\documentclass[letterpaper, 10pt, journal, twocolumn, final]{IEEEtran}
\usepackage{graphicx}
\usepackage{tikz}
\usepackage{amssymb,amsmath,cases}
\usepackage{algorithm}
\usepackage{algpseudocode}
\usepackage{algorithmicx}
\usepackage{amsfonts}
\usepackage[hidelinks]{hyperref}
\usepackage[normalem]{ulem}
\usepackage{cite}

\newcommand{\map}[3]{#1: #2 \rightarrow #3}
\newcommand{\N}{\mathbb{N}}
\newcommand{\R}{\mathbb{R}}
\newcommand{\set}{\{1, \dots, N\}}

\newcommand{\until}[1]{\{1,\ldots,#1\}}

\newcommand{\bL}{L}
\newcommand{\cL}{\mathcal{L}}
\newcommand{\tbL}{L}

\newcommand{\1}{\mathbf{1}}
\newcommand{\T}{^\top}
\newcommand{\inv}{^{-1}}
\newcommand{\norm}[1]{\left \|#1 \right \|}

\newcommand{\nablaF}{\nabla \mathbf{f}}
\newcommand{\nablatwoF}{\nabla^2 \mathbf{f}}

\newcommand{\tildex}{\tilde{x}} 
\newcommand{\tildez}{\tilde{z}} 
\newcommand{\tildey}{\tilde{y}} 
\newcommand{\tildeeta}{\tilde{\eta}} 
\newcommand{\tildepsi}{\tilde{\psi}} 

\newcommand{\xeq}{x_{\text{eq}}} 
\newcommand{\zeq}{z_{\text{eq}}}

\DeclareMathOperator{\col}{col}

\DeclareMathOperator{\diag}{diag}
\DeclareMathOperator{\blkdiag}{blkdiag}

\newtheorem{theorem}{Theorem}[section]
\newtheorem{proposition}[theorem]{Proposition}

\newtheorem{assumption}[theorem]{Assumption}
\newtheorem{remark}[theorem]{Remark} 
 \newtheorem{algo}[theorem]{Algorithm}

\newcommand\oprocendsymbol{\hbox{$\square$}}
\newcommand\oprocend{\relax\ifmmode\else\unskip\hfill\fi\oprocendsymbol}

\graphicspath{{figs/}}

\newcommand{\xstar}{x^\star}

\newcommand{\cD}{\mathcal{D}}
\newcommand{\cW}{\mathcal{W}}
\newcommand{\cS}{\mathcal{S}}
\newcommand{\EE}{\mathcal{E}}
\newcommand{\GG}{\mathcal{G}}

\newcommand{\NN}{\mathcal{N}}

\newcommand{\rx}{\mathrm{x}}
\newcommand{\rz}{\mathrm{z}}
\newcommand{\rs}{\mathrm{s}}

\newcommand{\xkp}{\rx_{k+1}}
\newcommand{\xk}{\rx_{k}}
\newcommand{\skp}{\rs_{k+1}}
\newcommand{\sk}{\rs_{k}}
\newcommand{\zkp}{\rz_{k+1}}
\newcommand{\zk}{\rz_{k}}

\newcommand{\rr}{\bar{r}}

\newcommand{\xikp}{\rx_{i,k+1}}
\newcommand{\xik}{\rx_{i,k}}
\newcommand{\xjk}{\rx_{j,k}}
\newcommand{\yikp}{\rs_{i,k+1}}
\newcommand{\yik}{\rs_{i,k}}
\newcommand{\yjk}{\rs_{j,k}}
\newcommand{\zikp}{\rz_{i,k+1}}
\newcommand{\zik}{\rz_{i,k}}
\newcommand{\zjk}{\rz_{j,k}}

\newcommand{\vx}{v_x}
\newcommand{\vz}{v_z}
\newcommand{\vg}{v_\nabla}
\newcommand{\vxi}{v_{i,x}}
\newcommand{\vzi}{v_{i,z}}
\newcommand{\vgi}{v_{i,\nabla}}
\newcommand{\vxz}{v_{xz}}
\newcommand{\vxzg}{v_{xz\nabla}}

\newcommand{\btk}{\tilde{t}^k}
\newcommand{\btkp}{\tilde{t}^{k+1}}
\newcommand{\tki}{t^{k_i}_i}
\newcommand{\tkj}{t^{k_j}_j}
\newcommand{\tkpi}{t^{k_i+1}_i}

\newcommand{\hxi}{\hat{x}_i^{k}}
\newcommand{\hxj}{\hat{x}_j^{k}}
\newcommand{\hnfi}{\nabla f_i^{k}}
\newcommand{\hnfj}{\nabla f_j^{k}}
\newcommand{\hzi}{\hat{z}_{i}^{k}}
\newcommand{\hzj}{\hat{z}_{j}^{k}}

\newcommand{\thxi}{\hat{\tildex}_i^{k}}
\newcommand{\thxj}{\hat{\tildex}_j^{k}}
\newcommand{\thzi}{\hat{\tildez}_{i}^{k}}
\newcommand{\thzj}{\hat{\tildez}_{j}^{k}}

\newcommand{\hxk}{\hat{x}^k}
\newcommand{\hzk}{\hat{z}^k}
\newcommand{\hnfk}{\nablaF^k}

\newcommand{\enf}{e_\nabla}

\newcommand{\kron}{\otimes}
\newcommand{\tr}{\lambda}
\allowdisplaybreaks

\def\algo/{{\scshape Continuous Gradient Tracking}}
\def\algoPeriodic/{{\scshape Synchronous Triggered Gradient Tracking}}
\def\algoTriggered/{{\scshape Asynchronous Triggered Gradient Tracking}}

\def\CGTP/{{\scshape Synchronous Triggered Gradient Tracking}}
\def\CGTT/{{\scshape Asynchronous Triggered Gradient Tracking}}

\def\er/{Erd\H{o}s-R\'enyi}

\begin{document}

\title{Triggered Gradient Tracking \\ for Asynchronous Distributed Optimization}
\author{Guido Carnevale, Ivano Notarnicola, Lorenzo Marconi,
Giuseppe Notarstefano
\thanks{The authors are with the Department of Electrical, 
Electronic and Information Engineering, University of Bologna, Bologna, Italy,
\texttt{\{name.lastname\}@unibo.it}.
This result is part of a project that has received funding from the European 
Research Council (ERC) under the European Union's Horizon 2020 research 
and innovation programme (grant agreement No 638992 - OPT4SMART).
}}

\maketitle
\begin{abstract}
	This paper proposes \algoTriggered/, i.e., a distributed optimization algorithm to solve consensus optimization over networks with asynchronous communication. As a building block, we devise the continuous-time counterpart of the recently proposed (discrete-time) distributed gradient tracking called \algo/. By using a Lyapunov approach, we prove exponential stability of the equilibrium corresponding to agents' estimates being consensual to the optimal solution, with arbitrary initialization of the local estimates. Then, we propose two triggered versions of the algorithm. In the first one, the agents continuously integrate their local dynamics and exchange with neighbors their current local variables in a synchronous way. In \algoTriggered/, we propose a totally asynchronous scheme in which each agent sends to neighbors its current local variables based on a triggering condition that depends on a locally verifiable condition. The triggering protocol preserves the linear convergence of the algorithm and avoids the Zeno behavior, i.e., an infinite number of triggering events over a finite interval of time is excluded. By using the stability analysis of \algo/ as a preparatory result, we show exponential stability of the equilibrium point holds for both triggered algorithms and any estimate initialization. Finally, the simulations validate the effectiveness of the proposed methods on a data analytics problem, showing also improved performance in terms of inter-agent communication.
\end{abstract}

\section{Introduction}	

As the devices with computation and communication capabilities are growing all around us, several contexts ranging from machine learning to autonomous vehicles and cooperative robotics, need to control networks of communication and computing agents. Tasks must be performed in an independent and cooperative way, without resorting to a centralized unit. Often, these tasks can be mathematically posed as distributed optimization problems. For this reason, distributed optimization has become an active research area, see~\cite{nedic2018distributed,notarstefano2019distributed,yang2019survey} for recent overviews.
Some applications tasks can be performed relying on distributed continuous-time optimization schemes. In~\cite{gharesifard2014distributed}, a distributed continuous-time optimization algorithm is proposed to solve a cost-coupled convex problem over a weighted digraph. A constrained convex problem is solved in~\cite{liu2015second} for a network of agents having local, second-order dynamics. In~\cite{zeng2017distributed}, a nonsmooth convex optimization problem with local constraints is solved by using a novel distributed continuous-time projected algorithm. In~\cite{lin2017distributed}, a distributed optimization problem with general differentiable convex objective functions is studied for continuous-time multi-agent systems with single-integrator dynamics. Objective functions subjected to bounds, equality, and inequality constraints have been taken into account in~\cite{yang2017multi} and addressed with a proportional-integral protocol. A distributed optimization problem on directed graphs with nonconvex local objective functions is tackled in~\cite{li2018distributed} by using an adaptive continuous-time algorithm. In~\cite{hatanaka2018passivity}, a passivity-based approach is used to prove the convergence of distributed continuous-time schemes for both unconstrained and constrained scenarios, also in presence of communication delays. In~\cite{li2020input}, a continuous-time distributed optimization algorithm is shown to possess exponential convergence properties by decomposing it into a set of input feedforward passive systems that interact with each other. Paper~\cite{moradian2022distributed} proposes a continuous-time optimization algorithm designed by taking inspiration from the existing discrete-time algorithm known as Newton-Raphson method.
In fact, a part of works is recently trying to study the convergence properties of dynamic systems representing the continuous counterpart of existing iterative optimization algorithms. 
This line of research starts with the work~\cite{su2014differential}, where a second order differential equation is presented as the continuous counterpart of Nesterov's accelerated gradient method. The authors of~\cite{wibisono2016variational} propose a systematic way to go from continuous-time curves generated by a Lagrangian functional to a family of discrete-time accelerated algorithms. A connection between a Lyapunov approach and the so-called estimating sequence analysis (typically adopted for momentum-based algorithms) is proposed in~\cite{wilson2021lyapunov} to analyze accelerated optimization methods. In~\cite{shi2021understanding}, high resolution ordinary differential equations are introduced to better understand the continuous-time counterpart of Nesterov algorithm and heavy ball algorithm. Paper~\cite{diakonikolas2019approximate} studies the first order mirror descent algorithm by deriving ordinary differential equations from duality gaps. 

Distributed continuous-time schemes rely on continuous-time communication among the network agents. Being such a mechanism not implementable in practice, a research effort has been also devoted to study continuous-time distributed algorithms with discrete-time communication.
In this regard, the authors of~\cite{kia2015distributed} propose a continuous-time optimization algorithm for solving a cost-coupled problem in a distributed manner. Moreover, they also propose variants of the plain algorithm in which the communication between agents occurs both in a periodic and event-triggered fashion. In~\cite{liu2016event}, the same scenario is addressed with a continuous-time algorithm in which the agents can only provide quantized information to neighbors with an event-triggered protocol. In~\cite{deng2017distributed}, a continuous-time algorithm with event-triggered communication inspired from~\cite{kia2015distributed} is extended to reject external disturbances by relying on internal model concepts. Paper~\cite{kajiyama2018distributed} employs a subgradient method is employed with an event-triggered communication policy. In~\cite{yi2018distributed} a continuous-time distributed optimization algorithm with second order dynamics both with continuous communication and event-triggered communication between agents is proposed. The work~\cite{zhao2018distributed} considers a quadratic optimization problem considered and a continuous-time algorithm with event-triggered communication is proposed to solve it. 
Finally, nonconvex optimization problems are addressed in~\cite{adachi2021distributed} with an event-triggered implementation of the distributed gradient descent.

The main contribution of this paper is the development of three distributed algorithms for strongly convex consensus optimization problems over networks along with the theoretical proof of their linear
convergence rate\footnote{``Linear convergence to an optimum'' in optimization and ``Exponential stability of an equilibrium'' in control are two expressions that refer to closely related concepts. In our paper they are equivalent.}. Specifically, the development of two of them rely on alternative communication paradigms that allow for a practical implementation of continuous-time-based schemes on real devices.
The first algorithm is called \algo/ and has been designed as the continuous-time counterpart of the recently proposed (discrete-time) distributed gradient tracking.
In the second algorithm, agents send their own variables to neighbors ruled by a synchronous triggering condition. 
Whereas, the third algorithm implements an asynchronous communication protocol among the agents that trigger according to locally verifiable triggering conditions.
By relying on a Lyapunov approach, we first develop a system theoretical stability analysis to show that the equilibrium of a dynamical system equivalent to \algo/ is globally exponentially stable. The equilibrium point corresponds to agents' solution estimates being consensual and equal to the optimal solution of the original optimization problem. Then, this result is extended to handle the stability analysis of the two triggered algorithms. Indeed,  the stability of their equilibrium is addressed by recasting the triggered algorithms as perturbed versions of \algo/. Under suitable conditions on the triggering protocol, exponential stability of \algo/ can be proved with arbitrary initial estimates.
The asynchronous version of the algorithm uses an additional auxiliary variable and a specific triggering condition. These features exclude the Zeno behavior (i.e., an infinite number of events in a finite period of time) while preserving the linear convergence. As a side result, we provide a robust stability guarantee of our algorithms against inexactness computations and/or communications.

As mentioned above, the discrete-time method called gradient tracking is instrumental for our novel triggered algorithms. It has been obtained as an extension of~\cite{nedic2009distributed,nedic2010constrained}, in which the gradient method is combined with consensus, by introducing a ``tracking action". Such tracking action is based on the dynamic average consensus (see~\cite{zhu2010discrete,kia2019tutorial}) in order to let each agent estimate the gradient of a total cost function, which is only partially locally known. There exist several variants of the gradient tracking, see~\cite{varagnolo2016newton,dilorenzo2016next,nedic2017achieving,qu2018harnessing,xu2017convergence,xi2017addopt,xin2018linear,scutari2019distributed,pu2020push,carnevale2020distributed}. 
A control-based analysis of this algorithm has been proposed in~\cite{bin2019system}, where a suitable change of coordinates is considered and turns out to be fundamental also for the analysis performed in this paper. Finally, the approach in~\cite{bin2019system} has been exploited to design gradient tracking algorithms with sparse (non-necessarily diagonal) gains in~\cite{carnevale2020enhanced}.
The paper is organized as follows. The problem setup is given in Section~\ref{sec:setup}. In Section~\ref{sec:C-GT} the \algo/ is derived along with its convergence properties. In Section~\ref{sec:dicrete_communication} the triggered algorithms are proposed and analyzed. Section~\ref{sec:discussion} discusses robustness aspects, while numerical simulations are provided in Section~\ref{sec:numerical_simulations}. The proofs are all deferred to the Appendix.

\paragraph*{Notation}
The matrix $M \in \R^{n\times n}$ is said to be Hurwitz if its spectrum, denoted as $\sigma(M)$, lies in the open left-half complex plane.
The identity matrix in $\R^{m\times m}$ is $I_m$, while $0_m$ is the all-zero matrix in $\R^{m\times m}$. 
The vector of $N$ ones is denoted by $1_N$, while $\1$ denotes $1_N \otimes I_d$ with $\kron$ being the Kronecker product. 
The vertical concatenation of the column vectors $v_1$ and $v_2$ is $\col (v_1, v_2)$. 
We denote as $\diag (m_1,\dots,m_N) \in \R^{N \times N}$ the diagonal matrix with diagonal elements $m_1, \dots, m_N$. Similarly, $\blkdiag (M_1, \dots, M_N)$ is the block diagonal matrix whose $i$-th block is $M_i \in \R^{d_i\times d_i}$.

\section{Problem Set-up and Preliminaries}
\label{sec:setup}

\subsection{Optimization and Network Setup}
In this paper, we consider a network of $N$ agents that aim at solving optimization problems in the form
\begin{align} \label{eq:problem}
  \min_{x\in\R^d} \: & \: \sum_{i=1}^N f_i(x),
\end{align}
with the cost function $\map{f_i}{\R^d}{\R}$ known to agent $i$ only, for all $i\in\until{N}$.
Let the following hold.
\begin{assumption}\label{ass:convexity}
	For all $i$, the function $\map{f_i}{\R^d}{\R}$ is strongly
	convex with coefficient $\alpha>0$.\oprocend
\end{assumption}
\begin{assumption}\label{ass:lipschitz}
	For all $i$, the function $\map{f_i}{\R^d}{\R}$ has Lipschitz 
	continuous gradient with constant $L>0$.\oprocend
\end{assumption}

Notice that Assumption~\ref{ass:convexity} ensures that problem~\eqref{eq:problem} has a unique optimal solution, denoted by $\xstar \in\R^d$.

The $N$ agents cooperate by exchanging information among neighbors to solve problem~\eqref{eq:problem}.
Their interaction is modeled as a graph $\GG = (\until{N}, \EE)$, with $\EE \subset \until{N}^2$ being the edge set.  %
If an edge $(i,j)$ belongs to $\EE$, then agents $i$ and $j$ can exchange information, otherwise not.
The set of neighbors of an agent $i$ is $\NN_i := \{j \in \set \mid (i,j) \in \EE\}$ and includes $i$ itself, i.e., $\GG$ contains self-loops.
We associate to the graph $\GG$ a symmetric weighted adjacency matrix $\cW \in\R^{N\times N}$ whose entries match the graph, i.e., $[\cW]_{ij} >0$ whenever $(i,j)\in \EE$ and $[\cW]_{ij} =0$ otherwise. %
The weighted degree of an agent is defined as $d_i = \sum_{j\in \NN_i}[\cW]_{ij}$. %
Finally, we associate to $\GG$ the so-called Laplacian matrix defined as $\cL:=\mathcal{D} - \cW$, where $\mathcal{D} := \text{diag}(d_1,\ldots,d_N) \in \R^{N \times N}$.
Let the following assumption of the network holds.
\begin{assumption}\label{ass:network}
	$\GG$ is undirected and connected.   \oprocend
\end{assumption}

\subsection{Discrete Gradient Tracking}\label{sec:D-GT}

Let us recall the discrete-time gradient tracking, recently proposed in
the literature~\cite{varagnolo2016newton,dilorenzo2016next,nedic2017achieving,qu2018harnessing,xu2017convergence,xi2017addopt,xin2018linear,scutari2019distributed,pu2020push,carnevale2020distributed}.
At iteration $k \in \N$, each agent $i$ maintains a local estimate $\rx_{i,k} \in \R^d$ of the optimal solution of problem~\eqref{eq:problem} and an auxiliary state $\rs_{i,k} \in \R^d$. These states are iteratively updated by agent $i$ based on \emph{(i)} the current value of the its local cost function gradient $\nabla f_i(\rx_{i,k})$ and \emph{(ii)} the information received from its neighbors according to 
\begin{subequations}\label{eq:local_s_GT}
	\begin{align}
	\xikp\!&=\!\sum_{j \in \NN_i} [\cW]_{ij}\xjk-\gamma\, \yik
	\\
	\yikp\!&=\!\sum_{j \in \NN_i}\![\cW]_{ij}\yjk \!+\! \nabla f_i(\xikp) \!- \!\nabla f_i(\xik),
	\end{align}
\end{subequations}
where $\gamma > 0$ is a stepsize, while each
$[\cW]_{ij}$ is the $(i,j)$-th element of a symmetric doubly stochastic matrix
$\cW$, i.e., such that $[\cW]_{ij}\ge0$ and $\sum_{i=1}^N [\cW]_{ij} = \sum_{i=1}^N [\cW]_{ji} = 1$ for all $j=\until{N}$. The distributed algorithm~\eqref{eq:local_s_GT} can be written into an
aggregate form as
\begin{subequations}\label{eq:s_GT}
	\begin{align}
	\xkp &= W \xk-\gamma\, \sk
	\label{eq:s_GT:x}
	\\
	\skp &= W \sk + \nablaF(\xkp) - \nablaF(\xk).
	\label{eq:s_GT:s}
	\end{align}
\end{subequations}
where $\xk := \col(\rx_{1,k}, \dots, \rx_{N,k})$,
$\sk := \col(\rs_{1,k}, \dots, \rs_{N,k})$, and
$\nablaF(\xk) := \col(\nabla f_1 (\rx_{1,k}), \dots, \nabla f_N (\rx_{N,k}))$
collect all the local variables of the agents while we set
$W := \cW \otimes I_d$. It can be proved that under the initialization $\rs_0 = \nablaF(\rx_0)$, the
algorithm in~\eqref{eq:s_GT} converges to a point corresponding to a consensual optimal solution of problem~\eqref{eq:problem},
see~\cite{varagnolo2016newton,dilorenzo2016next,nedic2017achieving,qu2018harnessing,xu2017convergence,xi2017addopt,xin2018linear,scutari2019distributed,pu2020push,carnevale2020distributed}.

\section{\algo/} 
\label{sec:C-GT}

In this section, we introduce \algo/ which is a novel continuous-time, distributed algorithm to solve problem~\eqref{eq:problem}.
We first show how to derive it as the continuous-time counterpart of the discrete-time algorithm recalled in~\eqref{eq:s_GT}.
Then, the convergence result to the optimal solution of problem~\eqref{eq:problem} is provided.

\subsection{From Discrete to Continuous}
\label{sec:analogy}

Let us first introduce an alternative, causal formulation of the (discrete) gradient tracking.
Consider the nonlinear change of coordinates proposed in~\cite{bin2019system} that moves the gradient term at $k+1$ to the left-hand side in~\eqref{eq:s_GT}. Therefore, let $\zik := \yik - \nabla f_i(\xik)$ and rewrite~\eqref{eq:s_GT} as
\begin{subequations}\label{eq:local_GT}
	\begin{align*}
	\xikp &=\! \sum_{j \in \NN_i} [\cW]_{ij} \xjk-\gamma\, \zik - \gamma\, \nabla f_i(\xik)
	\\
	\zikp &=\! \sum_{j \in \NN_i} [\cW]_{ij} \zjk 
	+\!
	\sum_{j \in \NN_i} [\cW]_{ij} \nabla f_j(\xjk) \!-\! \nabla f_i(\xik).
	\end{align*}
\end{subequations}
The latter can be also written into an aggregate form as
\begin{subequations}\label{eq:GT}
	\begin{align}
	\xkp &= W\xk-\gamma\, \zk -\gamma\, \nablaF (\xk)
	\\
	\zkp &= W \zk - (I_{Nd} - W) \nablaF(\xk),
	\end{align}
\end{subequations}
where $\zk := \col(\rz_{1,k}, \dots, \rz_{N,k})$ while $\xk$, $\nablaF(\xk)$, $W$ and $\gamma$ are the same as in Section~\ref{sec:D-GT}.

\begin{figure}[!htpb]
	\centering
	\includegraphics[scale=1]{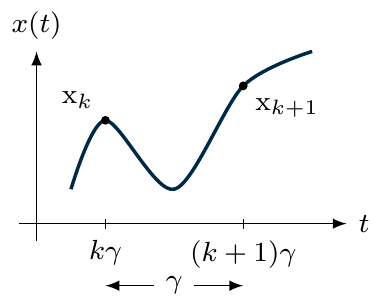}
	\caption{$\xk$ as sampled version of $x(t)$.}
	\label{fig:sampling}
\end{figure}

Following the arguments proposed in~\cite{su2014differential} for a centralized
optimization algorithm, we interpret the sequences $\xk$ and $\zk$ of~\eqref{eq:GT} as
sampled versions of two continuous-time signals $x(t)$ and $z(t)$.  These
signals are assumed to be smooth.
According to this interpretation, the stepsize $\gamma>0$ can be then seen as
the sampling time characterizing a discretization
procedure. Figure~\ref{fig:sampling} shows a graphical representation of the
discretization of the continuous signal $x(t)$ resulting in the discrete
sequence $\xk$.

Informally, we start from the intuition
\begin{align*}
	\xk \approx x(t)\Big|_{t=k\gamma}, 
	\qquad 
 	\zk \approx z(t)\Big|_{t=k\gamma},
\end{align*}
in which the discrete time $k$ is obtained by setting $k := t/\gamma$ with $t$ being the continuous time. For any fixed $t$, by choosing an arbitrarily small stepsize $\gamma$, we can consider the following approximations
\begin{align*}
x(t) &\approx \rx_{t/\gamma} = \xk,
\quad 
x(t + \gamma) \approx \rx_{(t + \gamma)/\gamma} = \xkp.
\end{align*}
The same clearly holds also for the sequence $\zk$. 
With these approximations in mind, we can write the following Taylor expansions
\begin{subequations}\label{eq:ansatz}
	\begin{align}
	\xkp \!&=\! x(t)\Big |_{t=(k+1)\gamma} \!\! =\!\! x(t)\Big |_{t=k\gamma} \!\! +\! \gamma \dot{x}(t)\Big |_{t=k\gamma} \!\!+\! o(\gamma)
	\\
	\zkp \!&=\! z(t)\Big |_{t=(k+1)\gamma} \!\! = \!\!z(t)\Big |_{t=k\gamma} \!\! +\! \gamma \dot{z}(t)\Big |_{t=k\gamma} \!\!+\! o(\gamma),
	\end{align}
\end{subequations}
where $o(\gamma)$ collects the higher order terms of the expansions.
As $\gamma$ goes to zero, the higher order terms in~\eqref{eq:ansatz} can be neglected, leading to
\begin{align*}
	\dot{x}(t) \!&=\! \tfrac{1}{\gamma} (\xkp-\xk) \!=\! \tfrac{W\!-\!I_{Nd}}{\gamma} x(t)\!-\!z(t)\!-\!\nablaF(x(t))
	\\
	\dot{z}(t) \!&=\! \tfrac{1}{\gamma} (\zkp-\zk) \!=\! \tfrac{W\!-\!I_{Nd}}{\gamma}z(t)\!-\!\tfrac{W\!-\!I_{Nd}}{\gamma}\nablaF(x(t)),
\end{align*}

which can be rewritten as
\begin{subequations}\label{eq:C-GT_sub_gamma}
	\begin{align}
	\dot{x}(t) &= -L_\gamma x(t) - z(t) - \nablaF(x(t))
	\\
	\dot{z}(t) &= -L_\gamma z(t) - L_\gamma\nablaF(x(t)),
	\end{align}
\end{subequations}
where $L_\gamma := (I_{Nd}-W)/\gamma$.
\subsection{Algorithm Definition and Convergence}

The Ordinary Differential Equation (ODE) in~\eqref{eq:C-GT_sub_gamma} involves matrices whose structure depends on the preceding derivation. However, one may consider any weighted Laplacian matrix $\cL$ associated to $\GG$. Therefore, we define the \algo/ as
\begin{align}\label{eq:C-GT}
\begin{bmatrix}
	\dot{x}(t)\\ \dot{z}(t)
\end{bmatrix} 
= 
\begin{bmatrix}
	-\bL& -I_{Nd}\\
	0& -\bL
\end{bmatrix}
\begin{bmatrix}
	x(t)\\ z(t)
\end{bmatrix}
-
\begin{bmatrix}
	I_{Nd} \\ \bL
\end{bmatrix} 
\nablaF(x(t)),
\end{align}
where $L \in \R^{Nd \times Nd}$ is given by $L := \cL \otimes I_d$.

It is useful to also provide a local view of algorithm~\eqref{eq:C-GT}, i.e., from the perspective of a generic agent $i$.
The $i$-th block-components of $x(t)$ and $z(t)$ corresponds, respectively, to the local states $x_i(t)$ and $z_i(t)$ of agent $i$.
The state $x_i(t)$ represents the local estimate at time $t$ of the optimal solution of problem~\eqref{eq:problem} while $z_i(t) \in \R^d$ is an auxiliary state. Set $\cW := \cD - \cL$ (with $\cD$ being the degree matrix of $\GG$) and let $w_{ij}$ be its $(i,j)$-th entry. 
Exploiting the sparsity in $\cW$, the $i$-th block-components of~\eqref{eq:C-GT} can be then written as
\begin{subequations}\label{eq:local_CGT}
\begin{align}
	\dot{x}_i (t)
	\!&=\!-\!\! 
	\sum_{j \in \NN_i} \! w_{ij} 
		\big(x_i(t)\!-\!x_j(t)\big) 
		\!-\! z_i(t)\! -\! \!\nabla f_i(x_i(t))
	\\
	\begin{split}
	\dot{z}_i(t) 
	\!&=\!-\!\!
	\sum_{j \in \NN_i} \! w_{ij} 
		\big(z_i(t) \!-\! z_j(t) \big) 
	\\
	&\hspace{1cm}
	-\!\! \sum_{j \in \NN_i} \! w_{ij} 
		\big( \nabla f_i(x_i(t)) - \nabla f_j(x_j(t)) \big).
	\end{split}
\end{align}
\end{subequations}

The following theorem establishes the convergence properties of \algo/.
\begin{theorem}\label{th:convergence}
  Consider the \algo/ distributed algorithm described by~\eqref{eq:C-GT}.  Let
  Assumptions~\ref{ass:convexity},~\ref{ass:lipschitz},~\ref{ass:network} hold
  and pick any $\col(x(0),z(0))$ such that $\1\T z(0) = 0$. Then, there exist $a_1 > 0$ and $a_2 >0$ such that
	\begin{align*}
		\norm{x_i(t) - \xstar} \leq a_1 \exp(-a_2t), \quad \forall i \in \set.\quad
		\oprocend
	\end{align*}
\end{theorem}
See Appendix~\ref{sec:proof_th1} for the proof. 

We underline that both $a_1$ and $a_2$ in Theorem~\ref{th:convergence} depend on \emph{(i)} the distance between the initial conditions and the system equilibrium and \emph{(ii)} the problem properties as, e.g., the network connectivity, the strong convexity parameter of the cost and the Lipschitz constants of the cost gradients. The same observation consistently applies to the subsequent results.

\begin{figure}[tpb]
	\centering
	\includegraphics[scale=1]{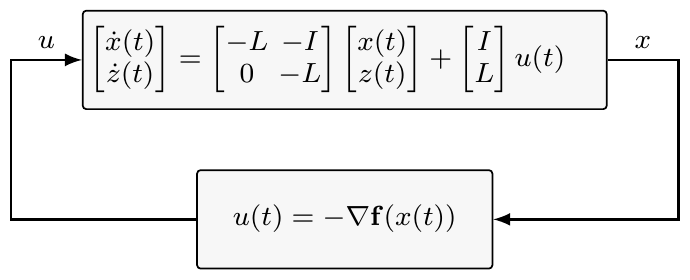}
	\caption{Block diagram representation of system~\eqref{eq:C-GT}.}
	\label{fig:block_diagram}
\end{figure}

We point out that the initialization $\1\T z(0) = 0$ can be obtained in a fully distributed way by simply setting each $z_i(0)=0$, for all $i \in \set$.
The proof of Theorem~\ref{th:convergence} is based on a Lyapunov analysis relying on the feedback structure of~\eqref{eq:C-GT} as represented in Figure~\ref{fig:block_diagram}.

Specifically, noticing that $\col(\1\xstar,-\nablaF(\1\xstar))$ is the unique
equilibrium point for system~\eqref{eq:C-GT} and by exploiting the initialization, we
perform a sequence of coordinate changes to obtain an equivalent, reduced system
formulation. The equivalent system is characterized by a (marginally stable)
linear part, associated to the consensus mechanism, which is perturbed by a
nonlinear (feedback) term depending on the gradient $\nablaF$.
We then prove the exponential stability of the equilibrium by designing a quadratic
Lyapunov function based on the linear part only and bounding the nonlinear
gradient term using the strong convexity and the Lipschitz continuity.

\begin{remark}
The expression of the discrete-time dynamics as in~\eqref{eq:GT} turns out to be crucial in the derivation of its continuous-time version. 
In fact, one can check that
\begin{align*}
	&\nablaF(\xkp) = \nablaF\left(x(t)\Big |_{t=(k+1)\gamma}\right) 
	\\
	&=  \nablaF\left( x(t)\Big |_{t=k\gamma} + \gamma \dot{x}(t)\Big |_{t=k\gamma} + o(\gamma)\right)
	\\
	&=  \nablaF\left( x(t)\Big |_{t=k\gamma}\right) 
	\!+\! \gamma \nablatwoF\left(x(t)\Big |_{t=k\gamma}\right)\dot{x}(t)\Big |_{t=k\gamma} \!+\! o(\gamma).
\end{align*}
Thus, the arguments presented in Section~\ref{sec:analogy} applied to the algorithm in its original coordinates~\eqref{eq:s_GT} results in
\begin{align*}
	\begin{bmatrix}
	\dot{x}(t)\\
	\dot{s}(t)
	\end{bmatrix} = \begin{bmatrix}
	-\bL_\gamma& -I_{Nd}\\
	- \nablatwoF(x(t))\tbL_\gamma& -(\tbL_\gamma + \nablatwoF(x(t)))
	\end{bmatrix}\begin{bmatrix}
	x(t)\\
	s(t)
	\end{bmatrix},
\end{align*}
where $s(t)$ is the continuous counterpart of $\sk$ while $\nablatwoF(x(t)) \!:=\!\blkdiag(\nabla^2 f_1(x_1(t)),\ldots,\nabla^2 f_N(x_N(t))) \in \R^{Nd \times Nd}$.
Notice that these coordinates involve the second-order matrix $\nablatwoF(\cdot)$ which usually requires a non-negligible computational complexity and, in certain applications, is not even known.
\oprocend
\end{remark}

\section{Triggered Gradient Tracking}
\label{sec:dicrete_communication}

The \algo/ would require communication among agents at all
$t\ge0$. Clearly, this prevents its practical implementation on real devices
that require time-slotted communication.
This issue is addressed next by proposing two alternative schemes in which inter-agent communication is triggered synchronously and asynchronously, respectively.

Let $\{\tki\}_{k_i \in \N}$, be the sequence of time instants at which agent $i$ sends its states $(x_i,z_i)$ and $\nabla f_i$ to its neighbors $j\in\NN_i$. 
Consistently, at time $\tkj$, agent $i$ receives the updated variables from its neighbor $j \in \NN_i$.
Let $\{\btk\}_{k\ge 0}$ be the ordered sequence of all the triggering times that occurred in the network. 
Then, given any $t \in [\btk,\btkp)$, let us introduce, for all $i \in\until{N}$, the shorthands
\begin{align}\label{eq:hat}
\begin{split}
	\hxi
	& := x_i(t) \Big|_{
		t= \inf \limits_{k_i \in \N}\left\{ \tki \ge \btk\right\}
		}
	\\
	\hzi
	&:=  z_i (t) \Big|_{
		t= \inf\limits_{k_i \in \N} \left\{ \tki \ge \btk\right\}
		}
	\\
	\hnfi
	& := \nabla f_i( x_i(t) ) \Big|_{
		t= \inf \limits_{k_i \in \N}\left\{ \tki \ge \btk\right\}
		}.
\end{split}
\end{align}
Quantities in~\eqref{eq:hat} represent the most updated values in the network within the considered time interval.
Under the described communication paradigm, we propose to modify the local dynamics in~\eqref{eq:local_CGT} as follows
\begin{subequations}\label{eq:periodic_local_CGT}
\begin{align}
	\dot{x}_i(t) \!&= \! - \!\!\! \sum_{j \in \NN_i} \! w_{ij}\! \left(\hxi - \hxj\right) - z_i(t) - \nabla f_i(x_i(t))
	\\
	\dot{z}_i(t) \!&= \! - \!\!\! \sum_{j \in \NN_i}\! w_{ij}\! \left(\hzi \! - \! \hzj\right) 
		\! - \!\!\! \sum_{j \in \NN_i} \! w_{ij}\! \left(\hnfi \! -  \!\hnfj\right)\!,
\end{align}
\end{subequations}
for all $t \in [\btk, \btkp)$. 
Within the $k$-th period, the variable $z_i$ behaves as an integrator.
As for the variable $x_i$, it is a local gradient flow compensated with an integral action $z_i$ and a constant consensus-error-like term. 
Agent $i$ does not use its own variables $x_i(t)$, $z_i(t)$, and $\nabla f_i(x_i(t))$ in the consensus mixing terms, but it rather uses their sampled version.
This fact allows one to preserve the theoretical consensus properties of the original scheme~\eqref{eq:local_CGT}.	

As one can expect, the specific rule to choose the triggering time $\tki$ will play a crucial role in the convergence properties of the resulting algorithms.

\subsection{\algoPeriodic/}
\label{sec:C-GT_synchronous}

We start by presenting \algoPeriodic/, obtained by imposing a \emph{synchronous
  communication} among agents.
Specifically, in this protocol each agent $i \in \set$ sends its local variables to its neighbors  at common instants of time chosen according to
\begin{align}\label{eq:periodic_law}
	\tkpi := \tki + \Delta,
\end{align}
for some common $\Delta > 0$ and with $t^{0_i} =t_0$ for all $i\in\set$. Intuitively,
the greater $\Delta$, the more inter-agent communication reduces. On the other
hand, the greater $\Delta$, the more the triggered algorithmic evolution moves
away from the behavior of \algo/. The next theorem gives theoretical guarantees
about the maximum admissible value for $\Delta$.
\begin{theorem}\label{th:periodic_convergence}
  Consider the algorithm in~\eqref{eq:periodic_local_CGT} with the synchronous
  communication protocol given by~\eqref{eq:periodic_law}.
  Let Assumptions~\ref{ass:convexity},~\ref{ass:lipschitz},~\ref{ass:network}
  hold and pick any $\col(x(0),z(0))$ such that $\1\T z(0) = 0$. Then, there
  exist $\Delta^\star > 0$, $a_3 >0$, and $a_4 > 0$ such that for any
    $\Delta \in (0, \Delta^\star)$, it holds
	\begin{align*}
		\norm{x_i(t) - \xstar} \leq a_3 \exp(-a_4t), \quad \forall i \in \set. \quad \oprocend%
	\end{align*}
\end{theorem}
See Appendix~\ref{sec:proof_periodic} for the proof.

For all $t \in [\btk,\btkp)$, the aggregate form of~\eqref{eq:periodic_local_CGT} reads as
\begin{align}
\label{eq:periodic_CGT_aggregate}
	\begin{bmatrix} 
		\dot{x}(t) \\ \dot{z}(t)
	\end{bmatrix} 
	& = 
	H
	\begin{bmatrix} 
		x(t) \\ z(t)
	\end{bmatrix}
	+
	B_1
	\nablaF(x(t))
	+
	B_2
	\begin{bmatrix} 
		\hxk \\ \hzk \\ \hnfk
	\end{bmatrix},
\end{align}
where $\hxk := \col(\hat{x}_1^{k},\dots,\hat{x}_N^{k})$, $\hzk := \col(\hat{z}_1^{k},\dots,\hat{z}_N^{k})$, $\hnfk := \col(\nabla f_1^{k},\dots,\nabla f_N^{k})$, and 
\begin{align*}
	H &:= 
	\begin{bmatrix}
		0& -I\\
		0& 0
	\end{bmatrix},
	 \quad
	B_1  := 
	\begin{bmatrix}
		-I \\
		0
	\end{bmatrix},
	 \quad
	B_2  := 
	\begin{bmatrix}
	-\bL& 0& 0\\
	   0& -\bL& -\bL
	\end{bmatrix}.
\end{align*}
The proof of Theorem~\ref{th:periodic_convergence} relies on a proper
reformulation of the dynamics~\eqref{eq:periodic_CGT_aggregate} as a perturbed
instance of the \algo/ system~\eqref{eq:C-GT}. In particular, we show that the
periodic triggering law~\eqref{eq:periodic_law} gives rise to a perturbation
term that vanishes at the equilibrium point (see, e.g.,~\cite[Chapter
9]{khalil2002nonlinear} for the notion of vanishing perturbation) and that can
be arbitrarily bounded through the parameter $\Delta$. Thus, we can consider the
same Lyapunov function $V$ used in the proof of Theorem~\ref{th:convergence} and
pick a sufficiently small $\Delta$ to show that the perturbation does not alter
the sign of the derivative of $V$ to conclude the proof.

\subsection{\algoTriggered/}
\label{sec:C-GT_asynchronous}

We now investigate the case in which the agents choose their triggering time
$\tki$ in a fully asynchronous way giving rise to an algorithm termed
\algoTriggered/. This scheme is motivated by the fact that the synchronous
communication executed according to~\eqref{eq:periodic_law} is rather
conservative with a consequent non-efficient usage of the available
resources. An asynchronous communication protocol allows agents to exchange
information only when really needed.
This requires a modification of the synchronous scheme. 
In particular, each agent has to check a local triggering condition and to maintain an additional auxiliary variable.
The latter is important to take into account the so-called Zeno behavior.
Specifically, an infinite number of triggerings over a finite interval of time must be avoided.
Indeed, for agent $i$, a triggering law suffers from the Zeno effect if 
\begin{align*}%
	\lim_{k_i\to\infty} \tki = \sum_{k_i=0}^\infty (\tkpi - \tki) = t^\infty_i,
\end{align*}
for some (finite) $t^\infty_i >0$ termed the Zeno time.

The local dynamics is again described by~\eqref{eq:periodic_local_CGT}. But, in order to perform communication only when needed, each agent $i$ chooses the next triggering time instant $\tkpi$ according to a locally verifiable condition. A possible choice for such a condition may be
\begin{align}
	\tkpi := \inf_{t > \tki} 
	\left \{
		\norm{e_i(t)} > \tr \norm{h_i(t)} \right \},\label{eq:triggering_without_xi}
\end{align}
with $e_i(t) := \col(x_i(t) - \hxi,z_i(t) - \hzi,\nabla f_i(x_i(t)) - \hnfi)$, $h_i(t) := z_i(t) + \nabla f_i(x_i(t))$, and $\tr > 0$ a constant to be properly specified later. The rationale for the triggering mechanism is to \emph{(i)} keep the triggered scheme close to the original dynamics~\eqref{eq:C-GT}, and \emph{(ii)} avoid the Zeno behavior. To this end, the right-hand side of the inequality within~\eqref{eq:triggering_without_xi} must be asymptotically vanishing when the algorithm approaches a steady state. This, in turn, gives rise to a vanishing quantity on the left term of the inequality. Indeed, looking also to the discrete-time version~\eqref{eq:local_s_GT}, the (local) quantity $z_i(t) + \nabla f_i(x_i(t))$ can be seen as a proxy for $\sum_{i=1}^N  \nabla f_i(x_i(t))$, i.e., a quantity that vanishes at a consensual optimal solution. 
However, $\norm{h_i(t)}$ vanishes not only when the algorithm approaches the equilibrium, but also if $(x_i(t),z_i(t)) \in \mathcal{S}_i := \{(x_i,z_i) \in \R^{2d} \mid z_i = -\nabla f_i(x_i)\}$, possibly giving rise to the Zeno behavior. 
Thus, in order to exclude this situation, the triggering condition~\eqref{eq:triggering_without_xi} is further modified as
\begin{align}\label{eq:triggering_law}
	\tkpi 
	:= 
	\inf_{t > \tki} 
	\Big \{
		\norm{e_i(t)} > \tr \norm{h_i(t)} 
	+ |\xi_i(t)| \Big \},
\end{align}
where $\xi_i \in \R$ is a local, auxiliary variable maintained by each agent $i$ evolving as 
\begin{align}\label{eq:xi}
	\dot{\xi}_i(t) = -\nu \xi_i(t),
\end{align}
where 
$\nu > 0$ is a parameter ruling the decay of $\xi_i(t)$.

As formally shown next, if the $\xi_i$ are initialized to nonzero values, then algorithm~\eqref{eq:periodic_local_CGT} with triggering law~\eqref{eq:triggering_law} does not incur in the Zeno behavior.
\begin{theorem}\label{th:triggered_convergence}
  Consider the algorithm described by~\eqref{eq:periodic_local_CGT} with the
  asynchronous communication protocol given by~\eqref{eq:triggering_law}.
	Let
  Assumptions~\ref{ass:convexity},~\ref{ass:lipschitz},~\ref{ass:network} hold
  and pick any $\col(x(0),z(0),\xi(0))$ such that $\1\T z(0) = 0$ and with $\xi(0) = \col(\xi_1(0),\ldots,\xi_N(0)) \ne 0$. Then, there exist $\tr^\star > 0$, $\nu^\star > 0$, $a_5 > 0$, and $a_6 >0$ such that for any $\tr \in (0, \tr^\star)$ in~\eqref{eq:triggering_law} and any $\nu > \nu^\star$, it holds
	\begin{align*}
		\norm{x_i(t) - \xstar} \leq a_5 \exp(-a_6t), \quad \forall i \in \set.%
	\end{align*}
	Moreover, system~\eqref{eq:periodic_local_CGT} does not exhibit the Zeno behavior.
	\oprocend
\end{theorem}
See Appendix~\ref{sec:proof_triggered} for the proof.

As for Theorem~\ref{th:periodic_convergence}, also this proof is based on a proper reformulation of
the aggregate form of \algoTriggered/ (which is still given
by~\eqref{eq:periodic_CGT_aggregate}) as a perturbed instance of the \algo/
dynamics~\eqref{eq:C-GT} with a vanishing perturbation.
For this asynchronous triggering law, \eqref{eq:triggering_law}, an upper bound on the perturbation magnitude is provided. It is proportional to \emph{(i)} the term $\lambda\norm{z(t) + \nablaF(x(t))}$, which, as already stated, represents a surrogate for the distance from the equilibrium point $\col(\1\xstar,-\nablaF(\1\xstar))$, and \emph{(ii)} to the exponentially decaying term $\norm{\xi}$.
Thus, considering a Lyapunov function derived from the one used in Theorem~\ref{th:convergence}, it is possible to show that, by picking suitable $\lambda$ and $\nu$, the perturbation does not affect the sign of the Lyapunov derivative.
A specific choice for $\tr^\star$ and $\nu^\star$ can be obtained by exploiting the bounds derived in the proof.

\section{Robustness Against Inexact Computation}
\label{sec:discussion}

Let us consider a more general scenario in which agents can access only inexact evaluations of their local state $(x_i,z_i)$ and/or of the local gradients $\nabla f_i$.
Let $\vgi (t)\in \R^d$ represents the mismatch between the exact value of $\nabla f_i(x_i(t))$ and the one available to agent $i$ for the local
updates. The presence of this mismatch may be due to several reasons as, e.g.,
quantization errors of the computing unit, measurement errors in the sensor
providing $\nabla f_i(x_i(t))$, or model uncertainties affecting the available
gradient.  Similarly, also mismatches affecting the states
$x_i$ and $z_i$ can be considered. Thus, we consistently introduce $\vxi (t) \in \R^d$ and
$\vzi(t) \in \R^d$ to model such uncertainties. This framework can be formalized
by writing
\begin{align}
\hspace*{-0.4cm}
  \begin{bmatrix}
    \dot{x} (t)
    \\
     \dot{z}(t)
  \end{bmatrix} 
  &=
  \begin{bmatrix}
  	-\bL  & -I\\
  	0  & -\bL
  \end{bmatrix}
	\begin{bmatrix}
    x(t) \\ z(t)
	\end{bmatrix}-
  \begin{bmatrix}
 		I \\ \bL
  \end{bmatrix} 
	\nablaF(x(t))
	\notag\\
	&\hspace{3.5cm} + \delta_1 B_2e(t)	+ B_3v(t),
	\label{eq:noise}
\end{align}
where $e(t) := \col(\hxk -x(t), \hzk - z(t), \hnfk - \nablaF(x(t)))$ collects (possible) mismatches due to discrete-time communication and $v(t):=\col(\vg(t),\vx(t),\vz(t))$ collects the mentioned local mismatches between the gradients, the solution estimates and the auxiliary variables, the matrix $B_2$ has the same meaning as in~\eqref{eq:periodic_CGT_aggregate}, and $B_3$ is defined as
\begin{align}
B_3 := \begin{bmatrix}
	-\bL& -I& -I\\
	0& -\bL& -\bL
	\end{bmatrix}.\label{eq:B3}
\end{align}
Finally, $\delta_1$ is equal to $0$ for \algo/ and equal to $1$ for both the \textsc{Synchronous} and \textsc{Asynchronous Triggered Gradient Tracking}.
Similarly, we denote as $\delta_2$ is equal to $0$ for both \algo/ and \textsc{Synchronous Triggered Gradient Tracking} and equal to $1$ for \textsc{Asynchronous Triggered Gradient Tracking}.

Next, the robustness of the algorithm in terms of input-to-state stability is studied.
\begin{proposition}\label{th:noise}
  Consider the algorithm described by~\eqref{eq:noise}.  Let
  Assumptions~\ref{ass:convexity},~\ref{ass:lipschitz},~\ref{ass:network} hold
  and pick any $\col(x(0),z(0))$ such that $\1\T z(0) = 0$.  Then, there exist a
  $\mathcal{KL}$ function $g_1(\cdot)$ and a $\mathcal{K}_\infty$ function
  $g_2(\cdot)$ such that for any $x(0) \in \R^{Nd}$ it holds
$
	\norm{x(t) - \1\xstar}\leq g_1(\norm{\chi(0)},t) + g_2(\norm{v(\cdot)}_\infty),
$
	with $\chi(0) :=\col(x(0)-\1\xstar,z(0) + \nablaF(\1\xstar),\delta_2\xi(0))$ and for any $v(\cdot) \in \mathcal{L}_\infty^{3Nd}$.\footnote{%
	See~\cite[Chapter 4]{khalil2002nonlinear} for the function classes' definitions.}
\oprocend
\end{proposition}
See Appendix~\ref{sec:proof_noise} for the proof.

Proposition~\ref{th:noise} guarantees that within the framework modeled by~\eqref{eq:noise}, the proposed algorithms behave as input-to-state stable systems. Therefore, in presence of mismatches on variables and gradients, the distance between the solution of problem~\eqref{eq:problem} and the computed estimates stay bounded according to the error magnitude.

\section{Numerical Simulations}
\label{sec:numerical_simulations}

We next present numerical simulations to confirm and support the theoretical findings.
The simulations are done using Matlab with its numerical solver ``ode45'' to integrate the \algo/.

We consider a network of agents that want to cooperatively solve a data analytics problem in which a linear classifier must be trained. %
Each agent $i$ is equipped with $m_i \in \mathbb{N}$ points $p_{i,1}, \dots, p_{i,m_i} \in \R^d$ with binary labels $l_{i,h} \in \{-1,1\}$ for all $h \in \{1,\dots,m_i\}$. 
We consider a logistic regression problem given by
\begin{align*}
	\min_{w,b}
	\sum_{i=1}^{N}
	\sum_{h=1}^{m_i} &\log\!\left(1\!+\!\exp ( -l_{i,h}(w\T p_{i,h} + b)) \right)\!+\!\tfrac{C(\norm{w}^2+b^2)}{2},
\end{align*}
where the optimization variables $w \in R^{d-1}$ and $b \in \R$ define the
  separating hyperplane, while $C > 0$ is the so-called regularization
parameter. Notice that the presence of the regularization makes the cost
function strongly convex.
In our simulations, we pick $d = 3$, $m_i = 10$ for all $i \in \{1,\dots, N \}$, and $C = 0.1$.

\subsection{\algo/}
\label{subsec:numerical_simulation_cgt}

\begin{figure}[!htpb]
	\centering
	 \includegraphics[scale=0.8]{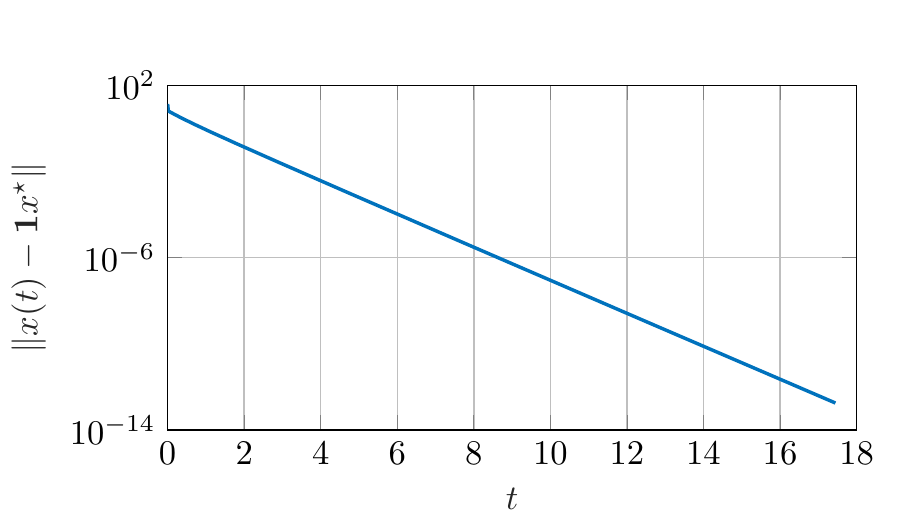}
	\caption{Evolution of the distance from the optimum of local estimates generated by \algo/.}
	\label{fig:CGT_error}
\end{figure}

In this subsection, the effectiveness of \algo/ is shown on a network of $N=50$ agents communicating according to an undirected and connected \er/ graph with parameter $0.4$. 
In Figure~\ref{fig:CGT_error} the convergence performances of \algo/ algorithm are shown.
Specifically, the distance of the local estimates $x(t) := \col(x_1(t), \ldots,x_N(t))$ from the optimum $\norm{x(t) - \1\xstar}$, converges to zero exponentially fast as expected from Theorem~\ref{th:convergence}.
\subsection{\scshape{Synchronous} and \scshape{Asynchronous Triggered Gradient Tracking}}

\begin{figure}[!htpb]
	\centering
	\includegraphics[scale=0.8]{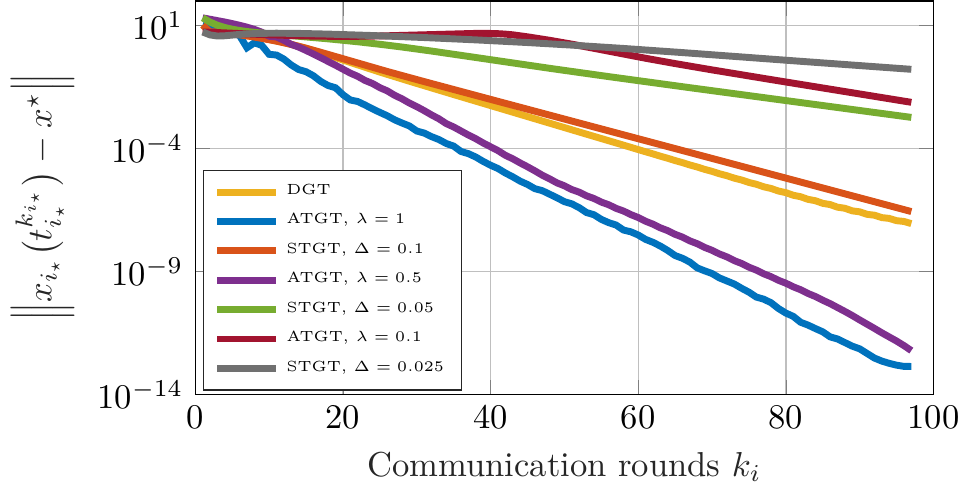}
	\caption{Comparison among \CGTT/ (ATGT), \CGTP/ (STGT) and the discrete gradient tracking (DGT) in terms of evolution of the optimality error.}
	\label{fig:triggered_error}
\end{figure}

In this subsection, the effectiveness of the triggered algorithms is shown for a
network of $N=10$ agents communicating according to an undirected and connected
\er/ graph with parameter $0.4$.
We tested \CGTP/ and \CGTT/ for different values of their key parameters $\Delta$ and $\lambda$, respectively. 
Moreover, we experimentally tuned the stepsize for the discrete graident tracking as $\gamma = 0.1$ in order to optimize its convergence rate. Finally, we set $\nu = 5$ for the dynamics of $\xi_i$ in~\eqref{eq:xi}. For the simulation of \CGTT/, the triggering condition (cf.~\eqref{eq:triggering_law}) is checked every $0.001$ seconds.
Figure~\ref{fig:triggered_error} compares the evolution of the optimality error obtained with different $\Delta$ and $\lambda$, for \CGTP/, \CGTT/, and the discrete gradient tracking. 
Specifically, the comparison is done in terms of communication rounds. 
The plot considers the performances of the most efficient agent, say $i_\star$, that performs the smallest number of neighboring communications in \CGTT/. As for the discrete gradient tracking, we denote, with a slight of abuse of notation, $x_{i_\star}(t_{i_\star}^{k_{i_\star}}) = \rx_{i_\star}^k$, with the sequence $\{\rx_{i_\star}^k\}_{k\ge 0}$ generated by~\eqref{eq:GT}.
As Figure~\ref{fig:triggered_error} clearly highlights, the communication rounds decrease as $\lambda$ increases. The same applies to $\Delta$. In particular, we underline that \CGTT/ results more efficient in finding the optimal solution with respect to both \CGTP/ and discrete gradient tracking.

\begin{figure}[!htpb]
	\centering
	\includegraphics[scale=0.8]{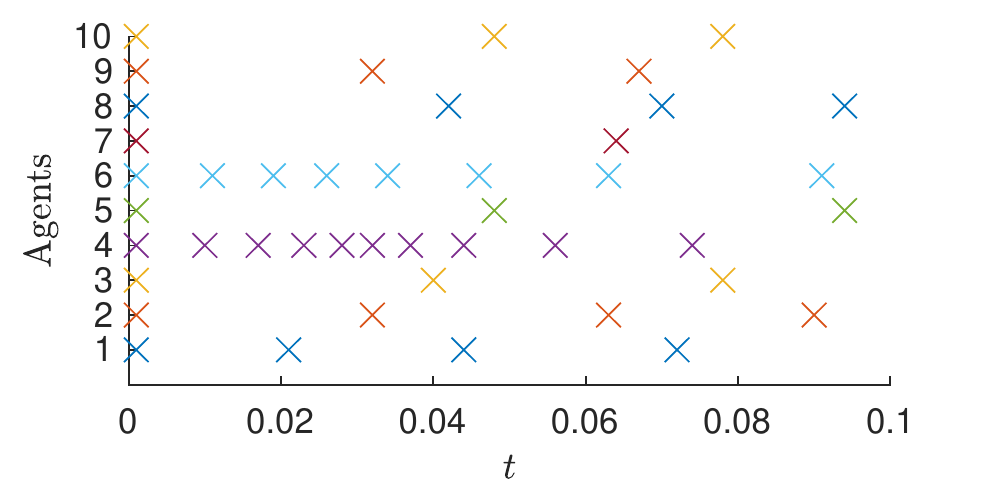}
	\caption{Occurrence of the triggering conditions in the \CGTT/.}
	\label{fig:communication_zoom}
\end{figure}

Finally, in Figure~\ref{fig:communication_zoom} each cross represents when the triggering condition occurred for each agent while running the \CGTT/ with $\tr = 0.1$. The plots demonstrate how event-triggered communication effectively reduces inter-agents communication.

\section{Conclusions}
\label{sec:conclusions}

In this paper, we addressed consensus optimization by proposing three novel
distributed optimization algorithms. First, \algo/ is derived as the
continuous-time counterpart of the existing discrete-time gradient tracking
algorithm. Then, by specifying proper inter-agent communication protocols, two triggered algorithms are
derived and analyzed, a synchronous and an asynchronous one, with an algorithmic structure inspired by the continuous version.  The convergence analysis
of all algorithms exploited a system-theoretical approach based on a suitably
defined quadratic Lyapunov function. %
The theoretical findings have been supported through numerical simulations.
\bibliographystyle{ieeetr}     

\begin{thebibliography}{10}

	\bibitem{nedic2018distributed}
	A.~Nedi{\'c} and J.~Liu, ``Distributed optimization for control,'' {\em Annual
		Review of Control, Robotics, and Autonomous Systems}, vol.~1, pp.~77--103,
		2018.
	
	\bibitem{notarstefano2019distributed}
	G.~Notarstefano, I.~Notarnicola, and A.~Camisa, ``Distributed optimization for
		smart cyber physical net-works,'' {\em Foundations and Trends in Systems and
		Control}, vol.~7, no.~3, pp.~253--383, 2019.
	
	\bibitem{yang2019survey}
	T.~Yang, X.~Yi, J.~Wu, Y.~Yuan, D.~Wu, Z.~Meng, Y.~Hong, H.~Wang, Z.~Lin, and
		K.~H. Johansson, ``A survey of distributed optimization,'' {\em Annual
		Reviews in Control}, vol.~47, pp.~278--305, 2019.
	
	\bibitem{gharesifard2014distributed}
	B.~Gharesifard and J.~Cort{\'e}s, ``Distributed continuous-time convex
		optimization on weight-balanced digraphs,'' {\em IEEE Transactions on
		Automatic Control}, vol.~59, no.~3, pp.~781--786, 2013.
	
	\bibitem{liu2015second}
	Q.~Liu and J.~Wang, ``A second-order multi-agent network for bound-constrained
		distributed optimization,'' {\em IEEE Transactions on Automatic Control},
		vol.~60, no.~12, pp.~3310--3315, 2015.
	
	\bibitem{zeng2017distributed}
	X.~Zeng, P.~Yi, and Y.~Hong, ``Distributed continuous-time algorithm for
		constrained convex optimizations via nonsmooth analysis approach,'' {\em IEEE
		Transactions on Automatic Control}, vol.~62, no.~10, pp.~5227--5233, 2017.
	
	\bibitem{lin2017distributed}
	P.~Lin, W.~Ren, and J.~A. Farrell, ``Distributed continuous-time optimization:
		nonuniform gradient gains, finite-time convergence, and convex constraint
		set,'' {\em IEEE Transactions on Automatic Control}, vol.~62, no.~5,
		pp.~2239--2253, 2017.
	
	\bibitem{yang2017multi}
	S.~Yang, Q.~Liu, and J.~Wang, ``A multi-agent system with a
		proportional-integral protocol for distributed constrained optimization,''
		{\em IEEE Transactions on Automatic Control}, vol.~62, no.~7, pp.~3461--3467,
		2017.
	
	\bibitem{li2018distributed}
	Z.~Li, Z.~Ding, J.~Sun, and Z.~Li, ``Distributed adaptive convex optimization
		on directed graphs via continuous-time algorithms,'' {\em IEEE Transactions
		on Automatic Control}, vol.~63, no.~5, pp.~1434--1441, 2018.
	
	\bibitem{hatanaka2018passivity}
	T.~Hatanaka, N.~Chopra, T.~Ishizaki, and N.~Li, ``Passivity-based distributed
		optimization with communication delays using pi consensus algorithm,'' {\em
		IEEE Transactions on Automatic Control}, vol.~63, no.~12, pp.~4421--4428,
		2018.
	
	\bibitem{li2020input}
	M.~Li, G.~Chesi, and Y.~Hong, ``Input-feedforward-passivity-based distributed
		optimization over jointly connected balanced digraphs,'' {\em IEEE
		Transactions on Automatic Control}, vol.~66, no.~9, pp.~4117--4131, 2021.
	
	\bibitem{moradian2022distributed}
	H.~Moradian and S.~S. Kia, ``A distributed continuous-time modified
		newton--raphson algorithm,'' {\em Automatica}, vol.~136, p.~109886, 2022.
	
	\bibitem{su2014differential}
	W.~Su, S.~Boyd, and E.~Candes, ``A differential equation for modeling
		nesterov’s accelerated gradient method: Theory and insights,'' in {\em
		Advances in Neural Information Processing Systems}, pp.~2510--2518, 2014.
	
	\bibitem{wibisono2016variational}
	A.~Wibisono, A.~C. Wilson, and M.~I. Jordan, ``A variational perspective on
		accelerated methods in optimization,'' {\em Proceedings of the National
		Academy of Sciences}, vol.~113, no.~47, pp.~E7351--E7358, 2016.
	
	\bibitem{wilson2021lyapunov}
	A.~C. Wilson, B.~Recht, and M.~I. Jordan, ``A lyapunov analysis of accelerated
		methods in optimization,'' {\em Journal of Machine Learning Research},
		vol.~22, no.~113, pp.~1--34, 2021.
	
	\bibitem{shi2021understanding}
	B.~Shi, S.~S. Du, M.~I. Jordan, and W.~J. Su, ``Understanding the acceleration
		phenomenon via high-resolution differential equations,'' {\em Mathematical
		Programming}, pp.~1--70, 2021.
	
	\bibitem{diakonikolas2019approximate}
	J.~Diakonikolas and L.~Orecchia, ``The approximate duality gap technique: A
		unified theory of first-order methods,'' {\em SIAM Journal on Optimization},
		vol.~29, no.~1, pp.~660--689, 2019.
	
	\bibitem{kia2015distributed}
	S.~S. Kia, J.~Cort{\'e}s, and S.~Mart{\'\i}nez, ``Distributed convex
		optimization via continuous-time coordination algorithms with discrete-time
		communication,'' {\em Automatica}, vol.~55, pp.~254--264, 2015.
	
	\bibitem{liu2016event}
	S.~Liu, L.~Xie, and D.~E. Quevedo, ``Event-triggered quantized
		communication-based distributed convex optimization,'' {\em IEEE Transactions
		on Control of Network Systems}, vol.~5, no.~1, pp.~167--178, 2016.
	
	\bibitem{deng2017distributed}
	Z.~Deng, X.~Wang, and Y.~Hong, ``Distributed optimisation design with triggers
		for disturbed continuous-time multi-agent systems,'' {\em IET Control Theory
		\& Applications}, vol.~11, no.~2, pp.~282--290, 2017.
	
	\bibitem{kajiyama2018distributed}
	Y.~Kajiyama, N.~Hayashi, and S.~Takai, ``Distributed subgradient method with
		edge-based event-triggered communication,'' {\em IEEE Transactions on
		Automatic Control}, vol.~63, no.~7, pp.~2248--2255, 2018.
	
	\bibitem{yi2018distributed}
	X.~Yi, L.~Yao, T.~Yang, J.~George, and K.~H. Johansson, ``Distributed
		optimization for second-order multi-agent systems with dynamic
		event-triggered communication,'' in {\em IEEE Conference on Decision and
		Control (CDC)}, pp.~3397--3402, 2018.
	
	\bibitem{zhao2018distributed}
	Z.~Zhao, G.~Chen, and M.~Dai, ``Distributed event-triggered scheme for a convex
		optimization problem in multi-agent systems,'' {\em Neurocomputing},
		vol.~284, pp.~90--98, 2018.
	
	\bibitem{adachi2021distributed}
	T.~Adachi, N.~Hayashi, and S.~Takai, ``Distributed gradient descent method with
		edge-based event-driven communication for non-convex optimization,'' {\em IET
		Control Theory \& Applications}, 2021.
	
	\bibitem{nedic2009distributed}
	A.~Nedi{\'c} and A.~Ozdaglar, ``Distributed subgradient methods for multi-agent
		optimization,'' {\em IEEE Transactions on Automatic Control}, vol.~54, no.~1,
		p.~48, 2009.
	
	\bibitem{nedic2010constrained}
	A.~Nedi{\'c}, A.~Ozdaglar, and P.~A. Parrilo, ``Constrained consensus and
		optimization in multi-agent networks,'' {\em IEEE Transactions on Automatic
		Control}, vol.~55, no.~4, pp.~922--938, 2010.
	
	\bibitem{zhu2010discrete}
	M.~Zhu and S.~Mart{\'\i}nez, ``Discrete-time dynamic average consensus,'' {\em
		Automatica}, vol.~46, no.~2, pp.~322--329, 2010.
	
	\bibitem{kia2019tutorial}
	S.~S. Kia, B.~Van~Scoy, J.~Cort{\'e}s, R.~A. Freeman, K.~M. Lynch, and
		S.~Mart{\'\i}nez, ``Tutorial on dynamic average consensus: The problem, its
		applications, and the algorithms,'' {\em IEEE Control Systems Magazine},
		vol.~39, no.~3, pp.~40--72, 2019.
	
	\bibitem{varagnolo2016newton}
	D.~Varagnolo, F.~Zanella, A.~Cenedese, G.~Pillonetto, and L.~Schenato,
		``{N}ewton-{R}aphson consensus for distributed convex optimization,'' {\em
		IEEE Transactions on Automatic Control}, vol.~61, no.~4, pp.~994--1009, 2016.
	
	\bibitem{dilorenzo2016next}
	P.~Di~Lorenzo and G.~Scutari, ``Next: In-network nonconvex optimization,'' {\em
		IEEE Transactions on Signal and Information Processing over Networks},
		vol.~2, no.~2, pp.~120--136, 2016.
	
	\bibitem{nedic2017achieving}
	A.~Nedi{\'c}, A.~Olshevsky, and W.~Shi, ``Achieving geometric convergence for
		distributed optimization over time-varying graphs,'' {\em SIAM Journal on
		Optimization}, vol.~27, no.~4, pp.~2597--2633, 2017.
	
	\bibitem{qu2018harnessing}
	G.~Qu and N.~Li, ``Harnessing smoothness to accelerate distributed
		optimization,'' {\em IEEE Transactions on Control of Network Systems},
		vol.~5, no.~3, pp.~1245--1260, 2018.
	
	\bibitem{xu2017convergence}
	J.~Xu, S.~Zhu, Y.~C. Soh, and L.~Xie, ``Convergence of asynchronous distributed
		gradient methods over stochastic networks,'' {\em IEEE Transactions on
		Automatic Control}, vol.~63, no.~2, pp.~434--448, 2017.
	
	\bibitem{xi2017addopt}
	C.~Xi, R.~Xin, and U.~A. Khan, ``{ADD-OPT}: Accelerated distributed directed
		optimization,'' {\em IEEE Transactions on Automatic Control}, vol.~63, no.~5,
		pp.~1329--1339, 2017.
	
	\bibitem{xin2018linear}
	R.~Xin and U.~A. Khan, ``A linear algorithm for optimization over directed
		graphs with geometric convergence,'' {\em IEEE Control Systems Letters},
		vol.~2, no.~3, pp.~315--320, 2018.
	
	\bibitem{scutari2019distributed}
	G.~Scutari and Y.~Sun, ``Distributed nonconvex constrained optimization over
		time-varying digraphs,'' {\em Mathematical Programming}, vol.~176, no.~1-2,
		pp.~497--544, 2019.
	
	\bibitem{pu2020push}
	S.~Pu, W.~Shi, J.~Xu, and A.~Nedi{\'c}, ``Push-pull gradient methods for
		distributed optimization in networks,'' {\em IEEE Transactions on Automatic
		Control}, vol.~66, no.~1, pp.~1--16, 2020.
	
	\bibitem{carnevale2020distributed}
	G.~Carnevale, F.~Farina, I.~Notarnicola, and G.~Notarstefano, ``Distributed
		online optimization via gradient tracking with adaptive momentum,'' {\em
		preprint arXiv:2009.01745}, 2020.
	
	\bibitem{bin2019system}
	M.~Bin, I.~Notarnicola, L.~Marconi, and G.~Notarstefano, ``A system theoretical
		perspective to gradient-tracking algorithms for distributed quadratic
		optimization,'' in {\em {IEEE} Conference on Decision and Control {(CDC)}},
		pp.~2994--2999, 2019.
	
	\bibitem{carnevale2020enhanced}
	G.~Carnevale, M.~Bin, I.~Notarnicola, L.~Marconi, and G.~Notarstefano,
		``Enhanced gradient tracking algorithms for distributed quadratic
		optimization via sparse gain design,'' {\em IFAC-PapersOnLine}, vol.~53,
		no.~2, pp.~2696--2701, 2020.
	
	\bibitem{khalil2002nonlinear}
	H.~K. Khalil, ``Nonlinear systems,'' {\em Upper Saddle River}, 2002.
	
	\bibitem{sontag2008input}
	E.~D. Sontag, ``Input to state stability: Basic concepts and results,'' in {\em
		Nonlinear and optimal control theory}, pp.~163--220, Springer, 2008.
	
	\end{thebibliography}

\appendix
\numberwithin{equation}{section}

\renewcommand{\thesection}{\Alph{section}}

\section{Proof of Theorem~\ref{th:convergence}}
\label{sec:proof_th1}
We first observe that, in light of Assumption~\ref{ass:lipschitz}, the ODE in~\eqref{eq:C-GT} is well posed for any $t\ge0$ and admits a
unique solution, see~\cite[Theorem~$3.2$]{khalil2002nonlinear}.
By inspecting system~\eqref{eq:C-GT}, we can assert that it has a unique
equilibrium point at
\begin{align*}
\begin{bmatrix}
	\xeq \\ \zeq
\end{bmatrix}
: = 
\begin{bmatrix}
	\1\xstar\\ -\nablaF(\1\xstar)
\end{bmatrix},
\end{align*}
which represents the situation in which the $N$ agents have a consensual solution
estimate equal to the optimal solution $\xstar$ of the optimization problem~\eqref{eq:problem}.
In order to use a Lyapunov approach, we put the system in error coordinates. Let
\begin{align}\label{eq:error_change}
\begin{bmatrix} x\\	z	\end{bmatrix} \longmapsto
	\begin{bmatrix}
		\tildex \\ \tildez
	\end{bmatrix} 
	:= 
	\begin{bmatrix} x\\	z	\end{bmatrix} 
	- 
	\begin{bmatrix} \xeq\\ \zeq \end{bmatrix}.
\end{align}
Then, system~\eqref{eq:C-GT} can be rewritten as
\begin{align}\label{eq:error_system}
\begin{bmatrix}
\dot{\tilde{x}}
\\
\dot{\tilde{z}}
\end{bmatrix} = \begin{bmatrix}
-\bL& -I\\
0& -\bL
\end{bmatrix}\begin{bmatrix}
\tilde{x}
\\
\tilde{z}
\end{bmatrix} + \begin{bmatrix}
I\\
\bL
\end{bmatrix}
u(\tildex),
\end{align}
where the role played by the ``input'' term $u(\tildex):= \nablaF(\1\xstar) - \nablaF(\tildex + \1\xstar)$ has been highlighted. Indeed, it can be interpreted as a nonlinear feedback of the output $\tildey = \tildex$ and suggests to introduce a further change of coordinates given by
\begin{align}\label{eq:normal_change}
	\begin{bmatrix} \tildex \\ \tildez \end{bmatrix} 
	\longmapsto
	\begin{bmatrix} 
		\tildey \\ \tildeeta
	\end{bmatrix} 
	:= 
	\underbrace{\begin{bmatrix}
	I & 0 \\
	\bL & -I
	\end{bmatrix}}_{T_1}
	\begin{bmatrix} \tildex \\ \tildez \end{bmatrix}.
\end{align}
Since $T_{1}$ is an involutory matrix (i.e., it coincides with its inverse), the change of 
coordinates~\eqref{eq:normal_change} transforms~\eqref{eq:error_system} in
\begin{align}\label{eq:normal_error_system}
	\begin{split}
	\begin{bmatrix}
	\dot{\tildey} \\
	\dot{\tildeeta}
	\end{bmatrix}
	& 
	= 
	\begin{bmatrix}
		-2\bL& I\\
		-\bL^2& 0 
	\end{bmatrix}
	\begin{bmatrix}
		\tildey \\ \tildeeta
	\end{bmatrix}  +
	\begin{bmatrix}
	I\\
	0
	\end{bmatrix}
	u(\tildey).
	\end{split}
\end{align}

Before studying the stability of the origin for~\eqref{eq:normal_error_system}, the effect of the initialization $\1\T z(0) = 0$ in the new coordinates $(\tildey,\tildeeta)$ is investigated.
We observe that the subspace 
\begin{align*}%
	\cS := \{(\tildey,\tildeeta) \mid  \1\T \tildeeta = 0 \}
\end{align*}
is invariant for~\eqref{eq:normal_error_system}. In light of~\eqref{eq:normal_change}, it holds
\begin{align*}
	0 
	& = \1\T \tildeeta 
	= \1\T (L \tildex - \tildez) 
	= \1\T z,
\end{align*}
where the last equality holds in light of~\eqref{eq:error_change} and since $\1\T \nablaF(\1\xstar)=0$ and $\1\T\bL = 0$ (cf.~Assumption~\ref{ass:network}).
Therefore the initialization of $z(0)$ guarantees that $\tildeeta(0) \in \cS$. 
Hence, we can perform a final change of coordinates to isolate the invariant state to further restrict the dynamics.
Let
\begin{align}\label{eq:mean_change_of_variables}
	\begin{bmatrix}
		\tildey \\ \tildeeta
	\end{bmatrix}
	\longmapsto 
	\begin{bmatrix}
		\tildey \\
		\tildepsi \\
		\tildeeta_{\text{avg}}
	\end{bmatrix}
	:= 
	T_2
	\begin{bmatrix}
		\tildey \\ \tildeeta
	\end{bmatrix}, 
\end{align}
in which 
\begin{align}\label{eq:T2}
	T_2 & := 
	\begin{bmatrix}
		T_{\tildey}\\T_{\tildeeta}
	\end{bmatrix}, \hspace{0.08cm}
	T_{\tildey} := 
	\begin{bmatrix}
		I & 0\\
		0 & R\T
	\end{bmatrix}, \hspace{0.08cm}
	T_{\tildeeta} := 
	\begin{bmatrix}
	0 & \tfrac{1}{\sqrt{N}} \1\T
	\end{bmatrix},
\end{align}
with $R \in \R^{Nd \times (N-1)d}$ such that $R\T R = I$, $R\T \1 = 0$ and $\norm{R} = 1$. The following useful relations holds true
\begin{align} \label{eq:RR}
	RR\T = I - \tfrac{1}{N}\1\1\T.
\end{align}
It is easy to check that $T_2\inv = T_2\T$, thus~\eqref{eq:normal_error_system} can be rewritten as
\begin{align}\label{eq:error_normalform}
	\begin{bmatrix}
		\dot{\tildey} \\
		\dot{\tildepsi} \\
		\dot{\tildeeta}_{\text{avg}}
	\end{bmatrix}
	\!&=\!\begin{bmatrix}
	-2\bL& R& \frac{\1}{\sqrt{N}}\\
	-R\T\bL^2& 0 &0\\
	0& 0& 0
	\end{bmatrix}\!
	\begin{bmatrix}
		\tildey \\
		\tildepsi \\
		\tildeeta_{\text{avg}}
	\end{bmatrix}\!+\!
	\begin{bmatrix}
	I\\
	0\\
	0
	\end{bmatrix}\!u(\tildey).
\end{align}
In light of the invariance of $\cS$, it holds $\tildeeta_{\text{avg}} (t) \equiv 0$. Then we can consider only $\zeta := \col(\tildey,\tildepsi) \in \R^n$, with $n := (2N-1)d$. The dynamics~\eqref{eq:error_normalform} can be written as
\begin{subequations}
\label{eq:system_A_B}
\begin{align}
	\label{eq:system_A_B_state}
	\dot{\zeta} &= A\zeta + Bu(\tildey),
\end{align}
\end{subequations}
with 
\begin{align}\label{eq:ABC_def}
	A &:= 
	\begin{bmatrix}
		-2\bL& R\\
		-R\T\bL^2& 0 
	\end{bmatrix},
	\quad
	B := \begin{bmatrix} I\\ 0	\end{bmatrix}.
\end{align}

Next, consider a quadratic, candidate Lyapunov function $V: \R^{n} \to \R$ given by
\begin{equation}\label{eq:lyapunov}
	V(\zeta) := \zeta\T P\zeta, 
\end{equation}
with $P \in \R^{n \times n}$ such that $P = P\T > 0$ and arranged in blocks as
\begin{equation}\label{eq:block_structure_P}
	P := 
	\begin{bmatrix}
		P_1& P_2\\
		P_2\T& P_3
	\end{bmatrix},
\end{equation}
where $P_1 \in \R^{Nd \times Nd}$, $P_2 \in \R^{Nd \times d}$, 
and $P_3 \in \R^{(Nd - d) \times (Nd - d)}$.
Next, it is shown how to choose $P$ in order to prove global exponential stability of the origin of~\eqref{eq:system_A_B}.
Let $m>0$ and set
\begin{align}\label{eq:P}
	P_1 &= m I, \quad P_2 &= - R,\quad P_3 &= m R\T (L^2)^\dagger R,
\end{align}
where $(\cdot)^\dagger$ denotes the Moore-Penrose pseudoinverse.
By the Schur complement lemma, $P>0$ imposes that $m$ must satisfy 
\begin{align}
\begin{cases}
	m I > 0\\
	mR\T (L^2)^\dagger R - \frac{1}{m}I > 0
\end{cases}
\Longrightarrow
m > \tfrac{1}{\sqrt{\min\{\sigma(R\T (\bL^2)^\dagger R)\}}}.
\label{eq:constraint_m1_P}
\end{align}

The time-derivative of $V$ along trajectories of~\eqref{eq:system_A_B} is
\begin{equation}\label{eq:dotV}
	\dot{V}(\zeta) = \zeta\T \underbrace{ (A\T P + PA) }_{-Q}\zeta + 2\zeta\T PBu.
\end{equation}
The choices~\eqref{eq:P} yield to
\begin{align}
	Q & =
	\begin{bmatrix}
		4 m \bL -2 \bL^2& 2\bL R\\
		2 R\T L& 2 I
	\end{bmatrix},
\quad
PB &= \begin{bmatrix}mI\\-R\T\end{bmatrix}.\label{eq:PB}
\end{align}
We separately study the quadratic term $-\zeta\T Q \zeta$ and the cross term $2\zeta\T PB u$ as a function of $m$ to show that $\dot{V}(\zeta)$ can be made negative definite for a sufficiently large $m$.
As for the first term in~\eqref{eq:dotV}, we observe that $Q$ is a solution to a Lyapunov equation associated to a marginally stable matrix. Therefore, it can only be positive semidefinite. Indeed, the upper-left block within the expression~\eqref{eq:PB} has the kernel spanned by $\1$ for any choice of $m$. %
By the Schur complement lemma, imposing $Q \ge 0$ is equivalent to
\begin{align}\label{eq:constraint_m1_Q_1}
 4m \bL -2 L^2 - 2\bL RR\T \bL \ge 0.
\end{align}
In light of~\eqref{eq:RR} and since $\bL\1 = 0$, condition~\eqref{eq:constraint_m1_Q_1} reduces to 
	$4m \bL -4 \bL^2 \ge 0$.
Since $\bL$ and $\bL^2$ have the same kernel, the latter condition is fulfilled by any $m$ such that
\begin{align}
	m \ge \dfrac{\max\{\sigma(\bL^2)\}}{\min\{\sigma(\bL) \setminus\{0\}\}}.\label{eq:constraint_m1_Q}
\end{align}
Moreover, $\bL$ positive semidefinite,
condition~\eqref{eq:constraint_m1_Q} can be satisfied with $m > 0$. 
Next, the second term in~\eqref{eq:dotV} is considered to show $\dot{V}<0$.  
In light of~\eqref{eq:PB}, it holds
\begin{align}
2\zeta\T PB u &= 2m\tildey \T u - 2\tildepsi \T R\T u\notag
\\
&\stackrel{(a)} \leq -2m\alpha \norm{\tildey}^2  - 2\tildepsi\T R\T u,\label{eq:cross_term_1}
\end{align}
where in $(a)$ we use the strong convexity of the cost functions (cf.~Assumption~\ref{ass:convexity}). 
Using the Cauchy-Schwarz inequality, condition~\eqref{eq:cross_term_1} can be manipulated as
\begin{align}
2\zeta\T PB u &\leq -2m\alpha \norm{\tildey}^2 + 2\norm{R}\norm{\tildepsi}\norm{u}\notag
\\
&\stackrel{(a)} \leq -2m\alpha\norm{\tildey}^2+ 2\beta \norm{\tildepsi}\norm{\tildey}\notag
\\
&\stackrel{(b)} \leq -2m\alpha \norm{\tildey}^2 + \frac{\beta }{\epsilon}\norm{\tildey}^2 + \beta \epsilon \norm{\tildepsi}^2\notag
\\
&\stackrel{(c)}= \zeta\T \underbrace{
	\begin{bmatrix}
 		\Big(-2m\alpha + \frac{\beta}{\epsilon} \Big) I& 0\\
 	0& \beta\epsilon  I
\end{bmatrix}
}_{Q_0} \zeta,
\label{eq:cross_term}
\end{align}
where in $(a)$ we use the Lipschitz continuity of the gradient of the cost functions (cf.~Assumption~\ref{ass:lipschitz}) and the fact that $\norm{R} = 1$, while in $(b)$ we use the Young's inequality with $\epsilon > 0$, and in $(c)$ we introduce the matrix $Q_0$. Indeed, we want to show that the zero eigenvalues of $Q$ can be moved inside the open left-half plane through $Q_0$. Thus, by plugging~\eqref{eq:cross_term} in~\eqref{eq:dotV}, it holds
\begin{align}
	\dot{V}(\zeta) \leq -\zeta\T\tilde{Q}\zeta,\label{eq:dotV2}
\end{align}
where $\tilde{Q} := Q - Q_0$, i.e., it holds
\begin{align}
\tilde{Q} \!&:=\!\!
 \begin{bmatrix}
4m \bL -2 \bL^2 + \Big(2m\alpha - \frac{\beta }{\epsilon} \Big) I& 2\bL R \\[0em]
2 R\T \bL& (2 - \beta\epsilon) I 
\end{bmatrix}.
\label{eq:Qtilde}
\end{align}
By the Schur complement lemma, $\tilde{Q} > 0$ is equivalent to
\begin{align*} %
	\begin{cases}
		 2 - \beta\epsilon > 0\\
		4m \bL - 2\left(\tfrac{2}{2-\beta\epsilon} +1\right) \bL^2 \!+\!  \tfrac{2m\alpha \epsilon -\beta }{\epsilon} I > 0,
	\end{cases}
\end{align*}
which is verified for every $\epsilon < \frac{2}{\beta}$ and
\begin{align} \label{eq:constraint_m1_tildeQ}
	m > \max
	\left\{
		\tfrac{\left(1+\frac{1}{2-\beta\epsilon}\right)\max\{\sigma(\bL^2)\}}{2\min\{\sigma(\bL)\setminus\{0\}\}}, 
	 \tfrac{\beta }{2\alpha\epsilon}
	 \right\}.
\end{align}
Therefore, we can conclude that $\dot{V}(\zeta) < -\min\{\sigma(\tilde{Q})\}\norm{\zeta}^2$ which implies that the origin is globally exponentially stable for system~\eqref{eq:system_A_B} (cf.~\cite[Theorem~4.10]{khalil2002nonlinear}).
Specifically, there exist $a_7, a_2 > 0$ such that
\begin{align}\label{eq:exponential}
\norm{\zeta(t)} \leq a_7\norm{\zeta(0)} \exp(-a_2t),
\end{align}
for any $\zeta(0) \in \R^n$. By noticing that 
$\norm{x_i(t) - \xstar} \leq \norm{x(t) - \1\xstar} = \norm{y(t)} \leq \norm{\zeta(t)},
$
the proof follows by~\eqref{eq:exponential} by setting $a_1 = a_7\norm{\zeta(0)}$.

\section{Proof of Theorem~\ref{th:periodic_convergence}}
\label{sec:proof_periodic}

The dynamics~\eqref{eq:periodic_local_CGT} of \algoPeriodic/ can be reformulated as a perturbed instance of the nominal dynamics of \algo/ described by~\eqref{eq:C-GT}. 
Clearly, the perturbation expresses the impact of the triggering mechanism on the algorithmic evolution. 
Thus, by adding and subtracting the term $B_2\col(x(t),z(t),\nablaF(x(t)))$ in the dynamics~\eqref{eq:periodic_CGT_aggregate}, we get
\begin{equation}\label{eq:periodic_CGT_error}
\begin{bmatrix}
	\dot{x} \\ \dot{z}
\end{bmatrix} = 
\begin{bmatrix} 
	-\bL& -I\\
	0& -\bL
\end{bmatrix}
\begin{bmatrix} 
	x \\ z
\end{bmatrix} 
+ 
\begin{bmatrix} 
	-I \\ -L
\end{bmatrix} \nablaF(x) 
+ B_2e,
\end{equation}
where $e$ has the same meaning as in~\eqref{eq:noise}.
By performing the same changes of coordinates defined in~\eqref{eq:error_change},~\eqref{eq:normal_change}, and~\eqref{eq:mean_change_of_variables}, the dynamics~\eqref{eq:periodic_CGT_error} can be equivalently reformulated as the following (restricted) dynamics
\begin{align}\label{eq:system_with_communication_error}
	\dot{\zeta} = A \zeta+ B u + E e_{\zeta,\nabla},
\end{align}
where the vectors $\zeta \in \R^n$, $u \in \R^{Nd}$ and the matrices $A \in \R^{n \times n}$ and $B \in \R^{n \times Nd}$ are as in~\eqref{eq:system_A_B}, while the quantities associated to the perturbation are
\begin{subequations}
	\begin{align}
	E & :=
	T_{\tildey}\T T_{1} B_2T_{1} T_{\tildey} 
	= 
	\begin{bmatrix}
		-\bL& 0& 0
		\\
		0 & -R\T\bL R& R\T\bL
	\end{bmatrix},\label{eq:definition_E}
	\\
	e_{\zeta,\nabla}
	& := 
	\col(
		\hat{\tildey} - \tildey,
		\hat{\tildepsi} - \tildepsi,
		\enf
	)
	\notag\\
	&:= 
	\col(
	\hat{\tildey} - \tildey,
	\hat{\tildepsi} - \tildepsi,
	\hnfk - \nablaF(\tildey+\xstar))
	\label{eq:definition_e}
\end{align}
\end{subequations}
with $T_{1}$ and $T_{\tildey}$ defined in~\eqref{eq:normal_change} and~\eqref{eq:T2}, respectively. 
Let us consider a quadratic, candidate Lyapunov function $V(\zeta) = \zeta \T P\zeta$
as in~\eqref{eq:lyapunov}
with the blocks of $P$ set as in~\eqref{eq:P}. The time-derivative of $V$ along the trajectories of~\eqref{eq:periodic_CGT_error} satisfies
\begin{align}
\dot{V}(\zeta) &= \zeta\T (A\T P + PA)\zeta + 2\zeta\T PB u + 2\zeta\T PEe_{\zeta,\nabla}
\notag\\
&\leq -\zeta\T \tilde{Q}\zeta + 2\zeta\T PEe_{\zeta,\nabla},\label{eq:periodic_inequality1}
\end{align}
where $\tilde{Q}$ is as in~\eqref{eq:Qtilde} so that the inequality holds in light the previous proof of Theorem~\ref{th:convergence} (cf.~\eqref{eq:dotV2}). 
By using the Young's inequality with $\epsilon > 0$, we can further upper bound~\eqref{eq:periodic_inequality1} as
\begin{align}
\dot{V}(\zeta)
& \leq -\zeta\T \tilde{Q}\zeta + \epsilon \zeta\T PP \zeta + \tfrac{1}{\epsilon}e_{\zeta,\nabla}\T E\T Ee_{\zeta,\nabla}
\notag\\
&\stackrel{(a)} = -\zeta\T \Big(\tilde{Q} - \epsilon P^2 \Big)\zeta + \tfrac{1}{\epsilon}e_{\zeta,\nabla}\T E\T Ee_{\zeta,\nabla},\label{eq:periodic_inequality2}
\end{align}
where in $(a)$ the terms in $\zeta$ have been grouped.
In light of the sufficient condition in~\eqref{eq:constraint_m1_tildeQ} to get a positive definite $\tilde{Q}$, we can always take $\epsilon$ such that
\begin{align*}
0 < \epsilon < \dfrac{\min\{\sigma(\tilde{Q})\}}{\max\{\sigma(P^2)\}}
\end{align*}
in order to impose also $\tilde{Q} - \epsilon P^2$ positive definite. Thus, by denoting as $q > 0$ the smallest eigenvalue of the matrix $\tilde{Q} - \epsilon P^2$ and by applying the Cauchy-Schwarz inequality to the quadratic term in $e_{\zeta,\nabla}$ of~\eqref{eq:periodic_inequality2}, we bound~\eqref{eq:periodic_inequality2} as
\begin{align}
\dot{V}(\zeta)& \leq -q\norm{\zeta}^2 + \tfrac{1}{\epsilon}\norm{E\T E}\norm{e_{\zeta,\nabla}}^2
\notag\\
& \stackrel{(a)}= -q\norm{\zeta}^2 + \tfrac{1}{\epsilon}\norm{E\T E}(\norm{e_\zeta}^2 + \norm{e_\nabla}^2)
\notag\\
& \stackrel{(b)}\leq -q\norm{\zeta}^2 + \tfrac{1}{\epsilon}\norm{E\T E}
	\Big( \norm{e_\zeta}^2 + \beta^2\norm{\hat{\tildey} - \tilde{y}}^2 \Big)
\notag\\
& \stackrel{(c)}\leq -q\norm{\zeta}^2 + \underbrace{ \tfrac{1}{\epsilon}\norm{E\T E}(1+\beta^2)}_{c_1}\norm{e_\zeta}^2,
\label{eq:vdot3}
\end{align}
where in $(a)$  we introduce $e_\zeta := \col(\hat{\tildey} - \tildey,
\hat{\tildepsi} - \tildepsi)$ to write $\norm{e_{\zeta,\nabla}}^2 = \norm{e_\zeta}^2 + \norm{e_\nabla}^2$, 
in $(b)$ we use the Lipschitz continuity of the gradients of the cost functions (cf.~Assumption~\ref{ass:lipschitz}) to bound $\norm{\enf}^2\leq \beta^2\| \hat{\tildey} - \tildey\|^2$, 
and in $(c)$ we rely on the fact that $\hat{\tildey} - \tilde{y}$ is a component of $e_\zeta$. 

The proof continues by deriving an upper bound for $\norm{e_\zeta}^2$ in~\eqref{eq:vdot3}. 
We start by defining
\begin{align}\label{eq:r_def}
	r := \dfrac{\norm{e_\zeta}}{\norm{\zeta}}.
\end{align}
Moreover, recall that in each interval $[\btk,\btkp)$, the error $e_\zeta$ is set to zero at $\btk$ and grows until $\btkp$ when it is reset again to zero. Hence, the goal is to establish a lower bound on the needed time for $r(t)$ to reach $\sqrt{q/c_1}$. By computing the time derivative of~\eqref{eq:r_def}, it follows 
\begin{align} 
\dot{r} &= \frac{e_\zeta\T \dot{e}_{\zeta}}{\norm{e_\zeta}\norm{\zeta}} - \frac{\norm{e_\zeta}\zeta\T\dot{\zeta}}{\norm{\zeta}^3}.\label{eq:rdot}
\end{align}
Using the Cauchy-Schwarz inequality, we bound $\dot{r}$ as
\begin{align}
\dot{r} 
& \leq 
\frac{\norm{e_\zeta} \norm{\dot{e}_{\zeta}}}{\norm{e_\zeta}\norm{\zeta}} 
+ 
\frac{\norm{e_\zeta}\norm{\zeta} \|\dot{\zeta}\| }{\norm{\zeta}^3}
\notag\\
&\stackrel{(a)} \leq 
\frac{ \|\dot{\zeta}\| }{\norm{\zeta}} 
+ 
\frac{\norm{e_\zeta} \|\dot{\zeta}\| }{\norm{\zeta}^2}\stackrel{(b)}{=}(1+r)\frac{ \|\dot{\zeta}\| }{\norm{\zeta}}
\label{eq:dotr1}
\end{align}
where in $(a)$ we use the identity $\dot{e}_{\zeta} = -\dot{\zeta}$ while 
in $(b)$ we exploit the definition of $r$ in~\eqref{eq:r_def}. Then, in light of the 
dynamics of $\zeta$
in~\eqref{eq:system_with_communication_error}, it holds
\begin{align}
&\dot{r} \leq (1+r)\frac{\norm{A\zeta+Bu+Ee_{\zeta,\nabla}}}{\norm{\zeta}}
\notag\\
&\stackrel{(a)}{\leq}
(1+r) \frac{
\norm{A} \norm{\zeta} + \norm{u} + \norm{E} \norm{e_{\zeta}} + \norm{E} \norm{e_\nabla}
}{\norm{\zeta}},
\label{eq:dotr2}
\end{align}
where in $(a)$ we use the triangle %
and the Cauchy-Schwarz inequalities combined with $\norm{Bu} = \norm{u}$. 
Next, by using the Lipschitz continuity of the gradients of the cost functions (cf. Assumption~\ref{ass:lipschitz}), we have $\norm{u} \leq \beta\norm{\zeta}$ and $\norm{e_\nabla} \leq \beta\norm{e_\zeta}$. Thus,~\eqref{eq:dotr2} becomes
\begin{align}
\dot{r} & \leq(1+r)\frac{(\norm{A}+\beta)\norm{\zeta} + (1+\beta)\norm{E}\norm{e_\zeta}}{\norm{\zeta}}
\notag\\
&\stackrel{(a)}{=} (1+r)\frac{\beta\norm{\zeta}}{\norm{\zeta}} + (1+r)\frac{\norm{A}\norm{\zeta} + (1+\beta)\norm{E}\norm{e_\zeta}}{\norm{\zeta}}
\notag\\
&\stackrel{(b)}{=} \beta(1+r) + (1+r)\frac{c_2\norm{\zeta} + c_2\norm{e_\zeta}}{\norm{\zeta}}
\notag\\
&\stackrel{(c)}{=} \beta(1+r) + c_2(1+r)^2,
\label{eq:dotr3}
\end{align}
where in $(a)$ we simply rearrange the terms, in $(b)$ we introduce $c_2 := \max\{\norm{A},(1+\beta)\norm{E}\}$, and in $(c)$ we use the definition of $r$ in~\eqref{eq:r_def}.

Using the Comparison Lemma (see~\cite[Lemma 3.4]{khalil2002nonlinear}) the bound~\eqref{eq:dotr3} translates in the following inequality
\begin{align}
	r(t,r(\btk)) \leq \rr(t,\rr(\btk)),\label{eq:r_s}
\end{align}
where $r(t,r(\btk))$ denotes the solution of~\eqref{eq:rdot} with initial condition at $t=t_k$ given by $r(t_k)$ 
while $\rr(t, \rr(\btk))$ denotes the solution of %
\begin{align}\label{eq:dots}
\dot{\rr}(t) = \beta(1+\rr(t)) + c_2(1+\rr(t))^2,
\end{align}
for some initial condition initial condition at $t=\btk$ given by $\rr(\btk)$ such that $r(\btk) \leq \rr(\btk)$.
Recalling that the protocol~\eqref{eq:periodic_law} imposes $r(\btk) = 0$ at the beginning of each time interval $[\btk,\btkp)$, then we select $\rr(\btk) = 0$. 
The solution of~\eqref{eq:dots} can be shown to be (cf.~\cite{kia2015distributed})
\begin{align}
 \rr(t,0)
	& = 
	\frac{(\beta + c_2)( \exp(\beta (t-\btk))  - 1)}{-c_2 \exp( \beta (t-\btk))
	+ 
	\beta + c_2 }.
	\label{eq:delta_solution}
\end{align}
Notice that $\rr(t,0)$ starts from $0$ at $t = \btk$ and monotonically increases within the interval 
$\Big[0, t_k + \ln\left(\frac{\beta+c_2}{c_2}\right)/\beta \Big)$. Thus, we can always find a triggering value $t = \Delta^\star>0$ such that $\rr(\Delta^\star,0) = \sqrt{ q/c_1}$. Hence, by choosing any $\Delta \in (0,\Delta^\star)$ in~\eqref{eq:periodic_law}, the inequality~\eqref{eq:r_s} ensures
\begin{align}
	|r(t)| = r(t) < \sqrt{\frac{q}{c_1}}, \label{eq:normr}
\end{align}
for all $t \in [\btk,\btkp)$, where the equality holds because $r$ is always positive, see its definition in~\eqref{eq:r_def}. 
With this result in mind, %
the inequality~\eqref{eq:vdot3} can be rewritten as
\begin{align*}
	\dot{V}(\zeta) \leq -\left(q - \frac{|r|^2}{c_1}\right)\norm{\zeta}^2, 
\end{align*}
 which allows us to use~\eqref{eq:normr} to conclude that the origin is globally exponentially stable for system~\eqref{eq:system_with_communication_error} (cf.~\cite[Th.~4.10]{khalil2002nonlinear}).
Specifically, there exist $a_4, a_8 >0$ such that
\begin{align}\label{eq:exponential_per}
\norm{\zeta(t)} \leq a_8\norm{\zeta(0)} \exp(-a_4t),
\end{align}
for any $\zeta(0) \in \R^n$. By noticing that
$
	\norm{x_i(t) - \xstar} \leq \norm{x(t) - \1\xstar} = \norm{y(t)} \leq \norm{\zeta(t)},
$
the proof follows by~\eqref{eq:exponential_per} by setting $a_3 = a_8\norm{\zeta(0)}$.

\section{Proof of Theorem~\ref{th:triggered_convergence} }
\label{sec:proof_triggered}

The proof of Theorem~\ref{th:triggered_convergence} traces the same initial steps in Section~\ref{sec:proof_periodic}. 
Specifically, we reformulate the \algoTriggered/ as a perturbed, extended version of \algo/ in which the perturbation is due to the event-triggered communication. By exploiting the steps leading to~\eqref{eq:system_with_communication_error}, the aggregate form of~\eqref{eq:periodic_local_CGT} and~\eqref{eq:xi} reads
\begin{subequations}\label{eq:triggered_perturbed_system}
\begin{align}
	\dot{\zeta} &= A\zeta + Bu + D e
	\\
	\dot{\xi} &= -\nu\xi,
\end{align}
\end{subequations}
where the vectors $\zeta \in \R^n$, $u \in \R^{Nd}$ and the matrices $A \in \R^{n \times n}$ and $B \in \R^{n \times Nd}$ are as in~\eqref{eq:system_A_B}, $e \in \R^{3Nd}$ has the same meaning as in~\eqref{eq:noise}, $\xi := \col(\xi_1,\dots,\xi_N) \in \R^N$, while the matrix $D$ is given by
\begin{equation}
	D :=T_{\tildey}\T T_{1} B_2
	= 
	\begin{bmatrix}
		-\bL& 0&0\\
		-R\T\bL^2& R\T\bL& R\T\bL
	\end{bmatrix},
\end{equation}
where the matrices $T_{1}$, $T_{\tildey}$ and $B_2$ are as in~\eqref{eq:normal_change},~\eqref{eq:T2}, and~\eqref{eq:periodic_CGT_aggregate}, respectively. We underline that the dynamics of $\zeta$ and $\xi$ are decoupled while both quantities affect the triggering law~\eqref{eq:triggering_law}.

Next, we show how to properly choose the value for $\nu$ in~\eqref{eq:xi} and for $\tr$ in the triggering law~\eqref{eq:triggering_law} to guarantee that the perturbation term $D e$ and the auxiliary variable $\xi$ do not alter the stability property associated the nominal system $\dot{\zeta} = A\zeta + Bu$ (cf.~Theorem~\ref{th:convergence}). 
To this end, an upper bound for $\norm{De}$, proportional to $\norm{\zeta}$ and $\norm{\xi}$, is derived. We start by using the Cauchy-Schwarz inequality to write $\norm{De} \leq \norm{D} \norm{e}\leq c_3\sum_{i=1}^N\norm{e_i}$, with $c_3 := \norm{D}$.
In light of the triggering law~\eqref{eq:triggering_law}, the latter inequality can be upper bounded as
\begin{align}
	\norm{De} & \leq \tr c_3 \sum_{i=1}^N\norm{z_i + \nabla f_i(x_i)} + c_3\sum_{i=1}^{N}|\xi_i|
	\notag\\
	&\stackrel{(a)}{\leq} \tr c_3 \sqrt{N}\norm{z + \nablaF(x)} + c_3\sqrt{N}\norm{\xi}
	 \notag\\
	 &\stackrel{(b)}{=} 
	\tr  c_4\norm{\tilde{z}\!+\!\nablaF(\tilde{x} + \1\xstar)\!-\! \nablaF(\1\xstar)}\!+\! c_4\norm{\xi} \notag\\
	&\stackrel{(c)}{\leq} 
	\tr  c_4\norm{\tilde{z}} + \lambda c_4\beta\norm{\tildex} + c_4\norm{\xi},\label{eq:triggering_inequality2}
\end{align}
where in $(a)$ we apply the basic algebraic relation $\sum_{i=1}^{N}\norm{\theta_i} \leq \sqrt{N}\norm{\theta}$ for a vector $\theta = \col(\theta_1,\dots,\theta_N)$, in $(b)$ we perform the change of coordinates given in~\eqref{eq:error_change} and introduce the constant $c_4 := c_3\sqrt{N}$, and in $(c)$ we use the triangle inequality and the Lipschitz continuity of the gradients of the cost functions (cf. Assumption~\ref{ass:lipschitz}). According to~\eqref{eq:normal_change} and~\eqref{eq:mean_change_of_variables}, it holds
\begin{align}
\begin{bmatrix}
\tilde{x}\\
\tilde{z}
\end{bmatrix} = T_1 T_2\T\begin{bmatrix}
\zeta\\
\tilde{\eta}_{\text{avg}}
\end{bmatrix} = T_1 T_2\T \begin{bmatrix}
\zeta\\
0
\end{bmatrix},\label{eq:zeta_xz}
\end{align}
where we use the fact that the initialization $z(0)$ leads to $\tilde{\eta}_{\text{avg}}(t) \equiv 0$. We rearrange the inequality~\eqref{eq:triggering_inequality2} to reconstruct the term $\norm{\col(\tildex,\tildez)}$ as
\begin{align}
\norm{De} & \leq \tr c_4\max\{1,\beta\}\sqrt{2}
\norm{\col(\tildex,\tildez)}
+ c_4\norm{\xi}
\notag\\
&\stackrel{(a)} \leq \tr c_5\norm{\zeta} + c_4\norm{\xi},\label{eq:inequality}
\end{align}
where in $(a)$ we combine~\eqref{eq:zeta_xz} with the Cauchy-Schwarz inequality and set $c_5 := c_4\max\{1,\beta\}\sqrt{2}\norm{T_1 T_{2}\T}$.
Given the linear bound in~\eqref{eq:inequality}, we can pursue a Lyapunov approach to conclude the global exponential stability of the origin. 
Let us consider a quadratic, candidate Lyapunov function $\tilde{V}(\zeta,\xi) = \zeta \T P\zeta + \tfrac{1}{2}\xi\T\xi$.
derived from the one considered in~\eqref{eq:lyapunov} with the blocks of $P$ set as in~\eqref{eq:P}. 
Using similar arguments leading to~\eqref{eq:vdot3}, the time-derivative of $\tilde{V}$ along trajectories of~\eqref{eq:triggered_perturbed_system} can be upper bounded as
\begin{equation}\label{eq:dotV_triggered}
	\dot{\tilde{V}}(\zeta,\xi) \leq -\tilde{q} \norm{\zeta}^2 + 2\zeta\T P D e - \nu\norm{\xi}^2.
\end{equation}
By using the Cauchy-Schwarz inequality, we can plug~\eqref{eq:inequality} in~\eqref{eq:dotV_triggered} to obtain
\begin{align}
	\dot{\tilde{V}}(\zeta,\xi)&\leq -\tilde{q} \norm{\zeta}^2 + 2\norm{\zeta}\norm{P} \norm{De}
	\notag
	\\
	&\le  -\tilde{q}_\tr\norm{\zeta}^2 +2c_4\norm{P}\norm{\zeta}\norm{\xi}- \nu\norm{\xi}^2,\label{eq:dotV_triggered2}
\end{align}
where we introduce $\tilde{q}_\tr := \left(\tilde{q} - 2\tr c_5\norm{P}\right)$. Then, for any $\tr < \frac{\tilde{q}}{2c_5}\norm{P} =: \tr^\star$, it holds $\tilde{q}_\tr > 0$. Setting $c_6 := c_4\norm{P}$, the inequality~\eqref{eq:dotV_triggered2} can be arranged in a matrix form as
\begin{align}
	\dot{\tilde{V}}(\zeta,\xi) \leq -
	\begin{bmatrix}
		\norm{\zeta}\\ \norm{\xi}
	\end{bmatrix}\T
		\underbrace{\begin{bmatrix}
		\tilde{q}_\tr& -c_6\\
		-c_6& \nu
		\end{bmatrix}}_{U}
	\begin{bmatrix}
		\norm{\zeta}\\ \norm{\xi}
	\end{bmatrix}.
	\label{eq:dotV_matrix}
\end{align}
Being $U \in \R^{2\times2}$ symmetric, by the Sylvester criterion $U > 0$ if and only if $\tilde{q}_\tr\nu > c_6^2$. Therefore, by taking any $\nu > \nu^\star :=\tfrac{c_6^2}{\tilde{q}_\tr} $, the matrix $U$ is positive definite. Thus the inequality~\eqref{eq:dotV_matrix} guarantees that %
the origin is globally exponentially stable for system~\eqref{eq:triggered_perturbed_system} (cf.~\cite[Lemma 4.10]{khalil2002nonlinear}).
Specifically, there exist $a_6, a_9 > 0$ such that
\begin{align}\label{eq:exponential_trigg}
\!\!\norm{\col(\zeta(t),\xi(t))}\!
\leq\!
\underbrace{a_9\!\norm{\col(\zeta(0),\xi(0))}}_{a_5} \exp(-a_6t),
\end{align}
for any $\col(\zeta(0),\xi(0)) \in \R^{N+n}$. By noticing that 
\begin{align*}
\norm{x_i(t) - \xstar} \leq \norm{x(t) - \1\xstar} \!=\! 
\norm{y(t)} 
\! \leq \! \norm{\col(\zeta(t),\xi(t))}\!,
\end{align*}
the proof of the first part of the theorem follows by~\eqref{eq:exponential_trigg}.

Next, we prove by contradction that~\eqref{eq:periodic_local_CGT} does not exhibit the Zeno behavior.
Suppose, without loss of generality, that an agent $i$ exhibits the Zeno behavior, namely 
\begin{align}\label{eq:zeno_lim}
	\lim_{k_i \to \infty} \tki = t_i^\infty.
\end{align}%
For any $t \ge 0$, we have
\begin{align}
	&\tfrac{d}{dt}\norm{e_i(t)}=\frac{e_i\T\dot{e}_i}{\norm{e_i(t)}}\stackrel{(a)}{\leq}\norm{\dot{e}_i(t)}
	\notag\\
	&\stackrel{(b)}{=}\norm{\col(
		\dot{x}_i(t),
		\dot{z}_i(t),
		\nabla^2 f_i(x_i(t))\dot{x}_i(t))}
	\notag\\
	&\stackrel{(c)}{=}\norm{
		\col(\dot{\tildex}_i(t),\dot{\tildez}_i(t),\nabla^2 f_i(\tildex_i(t) + \xstar)\dot{\tildex}_i(t))},\label{eq:dot_e_i}
\end{align}
where in $(a)$ we use the Cauchy-Schwarz inequality, in $(b)$ we use the definition of $e_i(t)$, and in $(c)$ we locally perform the change of variables given in~\eqref{eq:error_change}. 
Combining the latter change of variables with~\eqref{eq:periodic_local_CGT}, it holds
\begin{subequations}\label{eq:periodic_local_CGT_proof}
	\begin{align}
	\!\!\!\dot{\tildex}_i(t) \!\!
	&= 
	\!-\!\sum_{j \in \NN_i}\!\! w_{ij} (\thxi - \thxj) \!-\! \tildez_i(t) + u_i(\tildex_i(t))\label{eq:local_tildex}
	\\
	\!\!\!\dot{\tildez}_i(t) \!\!
	&= 
	\!-\!\sum_{j \in \NN_i}\!\!w_{ij} (\thzi\!-\!\thzj) 
	\!-\!\! \sum_{j \in \NN_i}\!\! w_{ij} (\hnfi\!-\!\hnfj),\label{eq:local_tildez}
	\end{align}
\end{subequations}
where we use $u_i(\tildex_i(t)) := (\nabla f_i(\tildex_i(t) + \xstar) - \nabla f_i(\xstar))$ and the local components of the shorthands given in~\eqref{eq:hat}. %
By~\eqref{eq:exponential_trigg}, the variables $\tildex_i(t)$ and $\tildez_i(t)$ are bounded for all $i\in \until{N}$ and $t \ge 0$.
Then, by defining $c_7:=\max_{i,t}\norm{\tildex_i(t)}$ and $c_8:=\max_{i,t}\norm{\tildez_i(t)}$,~\eqref{eq:local_tildex} and the triangle inequality can be combined to get 
\begin{align}
	\norm{\dot{\tildex}_i(t)} \! 
	\leq 
	\!\! \sum_{j \in \NN_i} \!\! w_{ij}2c_7 \!+\! c_8 \!+\! \norm{u_i(\tildex_i(t))}
	\!
	\stackrel{(a)}{\leq}
	\! 
	(2c_9\!+\!\beta) c_7 \!+\! c_8,
	\notag
\end{align}
where in $(a)$ we introduce $c_9 := \sum_{j \in \NN_i}w_{ij}$ and we use the Lipschitz continuity of the gradients of the cost functions (cf. Assumption~\ref{ass:lipschitz}). Using again the boundedness of the quantities, and by adding and subtracting $\nabla f_i(\xstar)$ within the second sum of~\eqref{eq:local_tildez}, it holds
\begin{align}\notag
	\norm{\dot{\tildez}_i(t)} \leq 2c_9(c_8 +\beta c_7).
\end{align} 
Moreover, the Lipschitz continuity of the gradients of the cost functions (cf. Assumption~\ref{ass:lipschitz}) also ensures that $\norm{\nabla^2f_i(v)} \leq \beta$, for all $v \in \R^d$ and all $i \in \until{N}$. By combining the latter with the two previous equations, the inequality~\eqref{eq:dot_e_i} can be upper bounded as
\begin{align}
	\tfrac{d}{dt}\norm{e_i(t)} \leq c_{10},\label{eq:bound_c_e}
\end{align}
with $c_{10} := (1+\beta)(2c_9+\beta)c_7 + c_8 + 2c_9(c_8 +\beta c_7)$.

Since the protocol~\eqref{eq:triggering_law} imposes $e_i(t) = 0$ at the beginning of each time interval $[\tki,\tkpi)$, then by also using~\eqref{eq:bound_c_e}, we can write
\begin{align}
	\!\!e_i(t) \!= \!e_i(\tki) + \int_{\tki}^t
	\!\dfrac{d\norm{e_i(\tau)}}{d\tau}d\tau \leq c_{10}(t-\tki).\label{eq:bound_e_i_final}
\end{align}
By~\eqref{eq:xi}, it holds $\xi_i(t) = \xi_i(0)\exp(-\nu t)$ for all $t \ge 0$. 
Thus, being $\lambda\norm{h_i(t)} \ge 0$ for any $t \ge 0$, the bound in~\eqref{eq:bound_e_i_final} imposes, as a necessary condition to satisfy the triggering in~\eqref{eq:triggering_law}, that
\begin{align}\label{eq:triggering_law_tilde}
	c_{10}(\tkpi - \tki) \ge |\xi_i(0)|\exp(-\nu \tkpi)
\end{align}
From~\eqref{eq:zeno_lim}, for all $\epsilon>0$ there exists $k_{i,\epsilon} \in \N$ such that
\begin{align}\label{eq:condition_limit}
	\tki \in [t_i^\infty - \epsilon, t_i^\infty], \quad \forall k_i \ge k_{i,\epsilon}.
\end{align}
Set
\begin{align}
	\epsilon := \frac{|\xi_i(0)|}{2c_{10}}\exp(-\nu t_i^\infty),
	\label{eq:epsilon}
\end{align}
and suppose that the $k_{i,\epsilon}$-th triggering time of agent $i$, namely $t_{i}^{k_{i,\epsilon}}$, has occurred. 
Let $t_i^{k_{i,\epsilon} + 1}$ be the next triggering time determined by~\eqref{eq:triggering_law}. Then, using the necessary condition~\eqref{eq:triggering_law_tilde} we can write
\begin{align}
	t_i^{k_{i,\epsilon} + 1} - t_i^{k_{i,\epsilon}} &\ge %
	\frac{|\xi_i(0)|}{c_{10}}\exp(-\nu t_i^{k_{i,\epsilon} + 1})
	\notag\\
	&\stackrel{(a)}{\geq} 
	\frac{|\xi_i(0)|}{c_{10}}\exp(-\nu t_i^\infty)
	\stackrel{(b)}{=} 
	2\epsilon,
	\label{eq:inequality_epsilon}
\end{align} 
where in $(a)$ %
we use $t_i^\infty \ge t_i^{k_{i,\epsilon} + 1}$, while in $(b)$ we use~\eqref{eq:epsilon}.
However~\eqref{eq:inequality_epsilon} implies
\begin{align*}
	t_i^{k_{i,\epsilon}} \leq t_i^{k_{i,\epsilon} + 1} - 2\epsilon \leq t_i^\infty - 2\epsilon,
\end{align*}
which contradicts~\eqref{eq:condition_limit} and concludes the proof.

\section{Proof of Proposition~\ref{th:noise}}
\label{sec:proof_noise}

The proof of Proposition~\ref{th:noise} traces the same initial steps in Section~\ref{sec:proof_periodic} and Section~\ref{sec:proof_triggered}.
Using the change of coordinates in~\eqref{eq:error_change},~\eqref{eq:normal_change},~\eqref{eq:mean_change_of_variables}, system~\eqref{eq:noise} can be recast as
\begin{align}\label{eq:changed_noise}
	\dot{\zeta} = A\zeta + B u + \delta_1 E e_{\zeta,\nabla} +  T_{\tildey}\T T_{1} B_3\vxzg,
\end{align}
with $\zeta \in \R^n$, $u \in \R^{Nd}$, where $A \in \R^{n \times n}$ and $B \in \R^{n \times Nd}$ are as in~\eqref{eq:system_A_B}, $E$ and $e_{\zeta,\nabla}$ are as in~\eqref{eq:definition_E} and~\eqref{eq:definition_e},
$B_3$, $T_{\tildey}$, and $T_1$ are as in~\eqref{eq:B3},~\eqref{eq:normal_change}, and~\eqref{eq:T2}, 
while $\vxzg := \col(\vx,\vz,\vg)$. We remark that $e_{\zeta,\nabla}$ changes according to the implemented communication protocol. Moreover, when \algoTriggered/ is considered, also dynamics~\eqref{eq:xi} has to be taken into account. However, when $\vg \equiv \vxz \equiv 0$, then $\vxzg \equiv 0$ and system~\eqref{eq:changed_noise} reduces to
\begin{align}\label{eq:equivalent_changed_noise}
\dot{\zeta} = A\zeta + B u + \delta_1 E e_{\zeta,\nabla}.
\end{align}
Theorems~\ref{th:convergence},~\ref{th:periodic_convergence}, and~\ref{th:triggered_convergence} ensure that the origin is globally exponentially stable for~\eqref{eq:equivalent_changed_noise} for both $\delta_1, \delta_2 \in \{0,1\}$ and for both communication protocols~\eqref{eq:periodic_law} and~\eqref{eq:triggering_law}.
In light of~\cite[Lemma 4.6]{khalil2002nonlinear}, this condition is sufficient to assert that system~\eqref{eq:changed_noise} is input-to-state stable and the proof follows (cf.~\cite[Section 2.9]{sontag2008input}).

\end{document}